\documentclass[a4paper, 14pt, fancybox]{article}
\raggedright
\usepackage{maple2e}
 
\DefineParaStyle{Maple Output}
\DefineParaStyle{Maple Output}
\DefineCharStyle{2D Math}
\DefineCharStyle{2D Output}

\usepackage[a4paper]{hyperref}

\usepackage{fancybox}
\usepackage{xy}
\usepackage[T1]{fontenc} 
\usepackage[latin1]{inputenc}
\usepackage{amsfonts,amssymb,amsmath}
\usepackage{graphics} 
\usepackage{multido, xcolor}

\newtheorem{theo}{\sc Theorem}[section]

\newtheorem{rema}[theo]{\sc Remark}

\newtheorem{conj}[theo]{\sc Conjecture}

\def\Im{{\hskip0.3mm\rm  Im\hskip0.5mm}}

\makeindex

\newcommand{\GL}{{\rm GL}}
\newcommand{\GSp}{{\rm GSp}}

\newcommand{\SL}{{\rm SL}}

\newcommand{\N}{{\mathbb N}}
\newcommand{\Q}{{\mathbb Q}}

\newcommand{\Z}{{\mathbb Z}}

\newcommand{\al}{\alpha}

\newcommand{\ep}{\varepsilon}

\font\teneusm=eusm10 \font\seveneusm=eusm7 
\font\fiveeusm=eusm5 
\newfam\eusmfam 
\textfont\eusmfam=\teneusm 
\scriptfont\eusmfam=\seveneusm 
\scriptscriptfont\eusmfam=\fiveeusm

\def\mat #1,#2,#3,#4,{\left({#1\atop #3}{#2\atop #4}\right)}
\def\bra#1,{{\left\lbrace {#1}\right\rbrace}}

\def \diag {{\rm diag}}

\def \Sp{{\rm Sp}}
\def\tm{\tilde m}

\def\l1{\langle}

\newcommand{\B}{\left(\begin{array}{cc}}
\newcommand{\E}{\end{array}\right)}

\def \fns{{${}^{*)}$}}

\newcommand{\comm}[1]
{\fns\marginpar{$\boxed
{\hskip-6pt
{\small {\sf 
\begin{tabular} {l}
 #1
\end{tabular}
}
}
}
$
}
}
\def \?  {\comm{check ?}}

\usepackage{euscript}
\let\scr=\EuScript

\let\mathcal=\scr           
\def\ang#1,{{\left\langle {#1}\right\rangle}} 
\newcommand{\ds}{\displaystyle}

\def\Im{{\hskip0.3mm\rm  Im\hskip0.5mm}}

\font\teneusm=eusm10 \font\seveneusm=eusm7 
\font\fiveeusm=eusm5 
\newfam\eusmfam 
\textfont\eusmfam=\teneusm 
\scriptfont\eusmfam=\seveneusm 
\scriptscriptfont\eusmfam=\fiveeusm

\font\tengothic=eufm10
\font\sevengothic=eufm7
\font\fivegothic=eufm5
\newfam\Gothic
\textfont\Gothic=\tengothic
\scriptfont\Gothic=\sevengothic
\scriptscriptfont\Gothic=\fivegothic

\def\mat #1,#2,#3,#4,{\left({#1\atop #3}{#2\atop #4}\right)}
\def\bra#1,{{\left\lbrace {#1}\right\rbrace}}

\def \diag {{\rm diag}}

\def \Sp{{\rm Sp}}

\def\l1{\langle}

\let\scr=\EuScript

\let\mathcal=\scr           

\def\dbC{{\db C}}

\def\vin{{ {\tiny \mid }  
\kern-7.29pt 
\bigcup }}

\def\ang#1,{{\left\langle {#1}\right\rangle}}

\normalsize

\newcommand{\CC}{\mathbb C}

\newcommand{\HH}{\mathbb H}

\newcommand{\QQ}{\mathbb Q}
\newcommand{\RR}{\mathbb R}

\newcounter{ncours}{\setcounter{ncours} {1}}

\def\diag {\mathop{\rm diag}\nolimits}
\def\SL {\mathop{\rm SL}\nolimits}
\def\Sp {\mathop{\rm Sp}\nolimits}

\input amssym.def
\input amssym.tex

\font\teneusm=eusm10 \font\seveneusm=eusm7 
\font\fiveeusm=eusm5 
\newfam\eusmfam 
\textfont\eusmfam=\teneusm 
\scriptfont\eusmfam=\seveneusm 
\scriptscriptfont\eusmfam=\fiveeusm 
\def\scr#1{{\fam\eusmfam\relax#1}}

\def\finishproclaim{\par\rm
    \ifdim\lastskip<\medskipamount\removelastskip
     \penalty55\medskip\fi}

\def\proofof#1:{\par\medskip
   \noindent{\it Proof of {\rm #1}}}

\def\Ref[#1]{\par\smallskip\hang\indent%
  \llap{\hbox to\parindent{[#1]\hfil\enspace}}%
     \ignorespaces}
\def\Item#1{\par\smallskip
  \hang\indent\llap{\hbox to\parindent
     {#1\hfill\enspace}}\ignorespaces}
\def\ItemItem#1{\par\indent\hangindent2\parindent
\hbox to\parindent{#1\hfill\enspace}\ignorespaces}

\def\arrowsim{\smash{\mathop{\longrightarrow}
 \limits^{\lower1.5pt \hbox{$\scriptstyle\sim$}}}}

 \def\diag{{\rm diag}}

 \def\SL{{\rm SL}}

\def\GL{{\rm GL}}

\def\dbC{{\Bbb C}}

\begin{document}

\title{
Explicit formulas for Hecke operators and Rankin's lemma in higher genus
}
\author{ Alexei Panchishkin, Kirill Vankov\\
{\tt http://www-fourier.ujf-grenoble.fr/\~{}panchish}
\\
\small e-mail : panchish$@$mozart.ujf-grenoble.fr, 
FAX:  33 (0) 4 76 51 44 78}

\date{}

\maketitle

 \begin{abstract}
We develop explicit formulas for Hecke operators of  higher genus
in terms of spherical coordinates.
Applications are given to summation of various generating series with coefficients in local Hecke algebra
and in a tensor product of  such algebras.
In particular, we formulate and prove Rankin's lemma in genus two. 
An application to a holomorphic lifting from $GSp_2 \times GSp_2$ to $GSp_4$ is given using Ikeda-Miyawaki constructions.


\end{abstract}
\tableofcontents

\section{Introduction: 
generating series for the  Hecke operators}

Let $p$ be a prime.
The Satake isomorphism \cite{Sa63} relates  $p$-local Hecke algebras of reductive groups over
$\Q$ to certain polynomial rings. 
Then one can  use a computer in order to find interesting identities between Hecke operators, 
between their eigenvalues, and  relations to  Fourier coefficients of modular forms of higher degree. 

The purpose of the present paper is to extend  Rankin's Lemma to the summation of Hecke series 
of higher genus using  symbolic computations.
We refer to \cite{Ma-Pa77}, where Rankin's Lemma was  used in the elliptic modular case
for  multiplicative and  additive convolutions of Dirichlet series.
That work was further developped in \cite{Pa87}, \cite{Pa02}, 
 see also \cite{Ma-Pa05}.


\

Recall that a classical method to produce $L$-functions for an algebraic group $G$ over $\Q$ 
uses the generating series 
$$
\sum_{n=1}^{\infty} \lambda_f(n)n^{-s}=\prod_{p \ {\rm primes}}\sum_{\delta=0}^{\infty} \lambda_f(p^\delta )p^{-\delta s},
$$
of the eigenvalues of Hecke operators
 on an automorphic form $f$ on $G$.
We study the generating series of Hecke operators $\mathbf{T}(n)$ for the symplectic group $Sp_g$, and $\lambda_f(n)=\lambda_f(\mathbf{T}(n))$.

Let  $\Gamma=\Sp_g(\Z)\subset \SL_{2g}(\Z)$ be the Siegel modular group of genus $g$, and  
 $[\mathbf{p}]_g=p\mathbf{I}_{2g}=\mathbf{T}(\underbrace{p, \cdots, p}_{2g})$  be the  scalar Hecke operator for $Sp_g$. 
According to Hecke and Shimura, 
\begin{align*}&
D_{p}(X)=\sum_{\delta=0}^{\infty} \mathbf{T}(p^\delta )X^\delta
\\ &
=\begin{cases}
\ds \frac 1 {1-\mathbf{T}(p)X+p[\mathbf{p}]_1X^2},  
\hskip 1.4cm  \mbox{ if } g=1 \\ 
\hskip 5cm
 \mbox{ (see \cite{Hecke}, and \cite{Shi71}, Theorem 3.21),}
\\ & \\
\ds \frac{1-p^2[\mathbf{p}]_2X^2}
 {
1-\mathbf{T}(p)X+\{p\mathbf{T}_1(p^2)+p(p^2+1)[\mathbf{p}]_2\}X^2
-p^3[\mathbf{p}]_2\mathbf{T}(p)X^3+p^6[\mathbf{p}]_2^2X^4
}&  \\ \hskip5cm\mbox{ if } g=2 \mbox{ (see \cite{Sh}, Theorem 2),}  
\end{cases}
\end{align*}
where 
$\mathbf{T}(p)$, $\mathbf{T}_i(p^2)$ ($i=1, \cdots, g$) are the $g+1$ generators of the corresponding Hecke ring over $\Z$
 for  the symplectic group $Sp_g$, in particular, $\mathbf{T}_g(p^2)=[\mathbf{p}]_g$.


In the present paper we study explicit formulas for Hecke operators in higher genus
in terms of spherical coordinates.
Applications are given to summation of various generating series with coefficients in local Hecke algebras.
The question of computing these series explicitely  was raised by Prof. S.Friedberg during 
first author's talk at the conference "Zeta Functions" (The Independent Moscow University, September 18-22, 2006).

In particular, we formulate and prove Rankin's lemma in genus two
for generating series with coefficients
 in a tensor product of  local Hecke algebras.

\section{Results}
\subsection{Preparation: a  formula for the total Hecke operator 
$\mathbf{T}(p^\delta)$  of genus 2}
We establish first the following  useful formula
 (in spherical variables $x_0, x_1, x_2$): 
\begin{align} \label{tpd}
{\Omega ^{(2)}_{x}}(\mathbf{T}&(p^{\delta }))  =
  p^{-1}\,\mathit{x}_0^{ \delta }
( p\,\mathit{x}_1^{(3+\delta)}\,\mathit{x}_2 
- p\,\mathit{x}_1^{(2+\delta)} 
- p\,\mathit{x}_1^{(3+\delta)}\,\mathit{x}_2^{(2+\delta)} 
+ p\,\mathit{x}_1^{(2+\delta)}\,\mathit{x}_2^{(3+\delta)} 
\\ & \quad \nonumber
- p\,\mathit{x}_1\,\mathit{x}_2^{(3+\delta)} 
+ p\,\mathit{x}_2^{(2+\delta)} 
+ p\,\mathit{x}_1 
- p\,\mathit{x}_2 
- \mathit{x}_1^{(2+\delta)}\,\mathit{x}_2^{2} 
+ \mathit{x}_1^{(1+\delta)}\,\mathit{x}_2 
\\ & \quad \nonumber
+ \mathit{x}_1^{(2+\delta)}\,\mathit{x}_2^{(1+\delta)} 
- \mathit{x}_1^{(1+\delta)}\,\mathit{x}_2^{(2+\delta)} 
+ \mathit{x}_1^{2}\,\mathit{x}_2^{(2+\delta)} 
- \mathit{x}_1\,\mathit{x}_2^{(1+\delta)} 
- \mathit{x}_1^{2}\,\mathit{x}_2 
+ \mathit{x}_1\,\mathit{x}_2^{2})/ 
\\ & \quad \nonumber
((1-\mathit{x}_1)\,(1-\mathit{x}_2)\,(1-\mathit{x}_1\,\mathit{x}_2)\,(\mathit{x}_1-\mathit{x}_2)) 
\end{align}
\begin{align*}
= & \nonumber
- p^{-1}\,\mathit{x}_0^{ \delta } 
((1-\mathit{x}_1\,\mathit{x}_2)\,(p\,\mathit{x}_1-\mathit{x}_2)\,\mathit{x}_1^{(\delta +1)} 
+(1-\mathit{x}_1\,\mathit{x}_2)\,(\mathit{x}_1-p\,\mathit{x}_2)\,\mathit{x}_2^{(\delta +1)} 
\\ & \quad \nonumber
-(1-p\,\mathit{x}_1\,\mathit{x}_2)\,(\mathit{x}_1-\mathit{x}_2)\,(\mathit{x}_1\,\mathit{x}_2)^{(\delta +1)}
-(p-\mathit{x}_1\,\mathit{x}_2)(\mathit{x}_1-\mathit{x}_2))/
\\ & \quad \nonumber
((1-\mathit{x}_1)\,(1-\mathit{x}_2)\,(1-\mathit{x}_1\,\mathit{x}_2)\,(\mathit{x}_1-\mathit{x}_2)). 
\end{align*}
\subsubsection*{Andrianov's generating series}
The expression (\ref{tpd}) comes from
the following  Andrianov's generating series  
\[
 \sum _{\delta =0}^{\infty } \,{\Omega ^{(2)}_{x}}(
\mathbf{T}(p^{\delta }))\,X^{\delta }=  {\displaystyle 
\frac {1-\frac{\mathit{x}_0^{2}\,\mathit{x}_1\,\mathit{x}_2}{p}\,X^{2} }
{(1-\mathit{x}_0\,X)\,(1-\mathit{x}_0\,\mathit{x}_1\,X )\,(1-\mathit{x}_0\,\mathit{x}_2\,X)\,
(1-\mathit{x}_0 \,\mathit{x}_1\,\mathit{x}_2\,X )}}
\]
after developing and a simplification using change of summation.

Note that the formula (\ref{tpd}) makes it possible to treat higher generating series of the following type
$$
 D_{p,m}(X)=\sum _{\delta =0}^{\infty } \,{\Omega ^{(2)}_{x}}(
\mathbf{T}(p^{m\delta }))\,X^{\delta } \ (m=2,3,\cdots)
$$
(in spherical variables $x_0, x_1, x_2$). 

\subsection{Rankin's generating series in genus 2} 
Let us use the spherical variables $x_0, x_1, x_2$ and $y_0, y_1, y_2$
for the Hecke operators.

Note that there are two types of convolutions:
the first one is  defined through the Fourier coefficients
(it was used by \cite{And-Kal78} for the analytic continuation of the standard $L$- function),
and the second one is defined through the eigenvalues of Hecke operators, 
and it is more suitable in order to treat the $L$-functions attached to tensor products of representations of the Langlands group.
However, a link  between the two types is known only for $g=1$.

In order to state a multiplicative analogue  of Rankin's lemma in genus two we need to write
the corrseponding formula for Hecke operator 
$\mathbf{T}(p^\delta)$  (in spherical variables $y_0, y_1, y_2$).

\begin{align*} 
{\Omega ^{(2)}_{y}}(\mathbf{T}&(p^{\delta }))  =
  p^{-1}\,\mathit{y}_0^{ \delta }
( p\,\mathit{y}_1^{(3+\delta)}\,\mathit{y}_2 
- p\,\mathit{y}_1^{(2+\delta)} 
- p\,\mathit{y}_1^{(3+\delta)}\,\mathit{y}_2^{(2+\delta)} 
+ p\,\mathit{y}_1^{(2+\delta)}\,\mathit{y}_2^{(3+\delta)} 
\\ & \quad 
- p\,\mathit{y}_1\,\mathit{y}_2^{(3+\delta)} 
+ p\,\mathit{y}_2^{(2+\delta)} 
+ p\,\mathit{y}_1 
- p\,\mathit{y}_2 
- \mathit{y}_1^{(2+\delta)}\,\mathit{y}_2^{2} 
+ \mathit{y}_1^{(1+\delta)}\,\mathit{y}_2 
\\ & \quad 
+ \mathit{y}_1^{(2+\delta)}\,\mathit{y}_2^{(1+\delta)} 
- \mathit{y}_1^{(1+\delta)}\,\mathit{y}_2^{(2+\delta)} 
+ \mathit{y}_1^{2}\,\mathit{y}_2^{(2+\delta)} 
- \mathit{y}_1\,\mathit{y}_2^{(1+\delta)} 
- \mathit{y}_1^{2}\,\mathit{y}_2 
+ \mathit{y}_1\,\mathit{y}_2^{2})/ 
\\ & \quad 
((1-\mathit{y}_1)\,(1-\mathit{y}_2)\,(1-\mathit{y}_1\,\mathit{y}_2)\,(\mathit{y}_1-\mathit{y}_2)) \,.
\end{align*}
Then we have that the product of the above polynomials is given by
\begin{align*}
{\Omega ^{(2)}_{x}}&(\mathbf{T}(p^{\delta })) \cdot{\Omega ^{(2)}_{y}}(\mathbf{T}(p^{\delta }))=
\\ & 
p^{-2}\,\mathit{x}_0^{\delta }\,\mathit{y}_0^{\delta }
( p\,\mathit{x}_1^{(3+\delta)}\,\mathit{x}_2 
- p\,\mathit{x}_1^{(2+\delta)} 
- p\,\mathit{x}_1^{(3+\delta)}\,\mathit{x}_2^{(2+\delta)} 
+ p\,\mathit{x}_1^{(2+\delta)}\,\mathit{x}_2^{(3+\delta)} 
\\ & \quad 
- p\,\mathit{x}_1\,\mathit{x}_2^{(3+\delta)} 
+ p\,\mathit{x}_2^{(2+\delta)} 
+ p\,\mathit{x}_1 
- p\,\mathit{x}_2 
- \mathit{x}_1^{(2+\delta)}\,\mathit{x}_2^{2} 
+ \mathit{x}_1^{(1+\delta)}\,\mathit{x}_2 
\\ & \quad 
+ \mathit{x}_1^{(2+\delta)}\,\mathit{x}_2^{(1+\delta)} 
- \mathit{x}_1^{(1+\delta)}\,\mathit{x}_2^{(2+\delta)} 
+ \mathit{x}_1^{2}\,\mathit{x}_2^{(2+\delta)} 
- \mathit{x}_1\,\mathit{x}_2^{(1+\delta)} 
- \mathit{x}_1^{2}\,\mathit{x}_2 
+ \mathit{x}_1\,\mathit{x}_2^{2})
\\ & \quad 
\times
( p\,\mathit{y}_1^{(3+\delta)}\,\mathit{y}_2 
- p\,\mathit{y}_1^{(2+\delta)} 
- p\,\mathit{y}_1^{(3+\delta)}\,\mathit{y}_2^{(2+\delta)} 
+ p\,\mathit{y}_1^{(2+\delta)}\,\mathit{y}_2^{(3+\delta)} 
\\ & \quad 
- p\,\mathit{y}_1\,\mathit{y}_2^{(3+\delta)} 
+ p\,\mathit{y}_2^{(2+\delta)} 
+ p\,\mathit{y}_1 
- p\,\mathit{y}_2 
- \mathit{y}_1^{(2+\delta)}\,\mathit{y}_2^{2} 
+ \mathit{y}_1^{(1+\delta)}\,\mathit{y}_2 
\\ & \quad 
+ \mathit{y}_1^{(2+\delta)}\,\mathit{y}_2^{(1+\delta)} 
- \mathit{y}_1^{(1+\delta)}\,\mathit{y}_2^{(2+\delta)} 
+ \mathit{y}_1^{2}\,\mathit{y}_2^{(2+\delta)} 
- \mathit{y}_1\,\mathit{y}_2^{(1+\delta)} 
- \mathit{y}_1^{2}\,\mathit{y}_2 
+ \mathit{y}_1\,\mathit{y}_2^{2})
/ 
\\ & \quad 
((1-\mathit{x}_1)\,(1-\mathit{x}_2)\,(1-\mathit{x}_1\,\mathit{x}_2)\,(\mathit{x}_1-\mathit{x}_2)
(1-\mathit{y}_1)\,(1-\mathit{y}_2)\,(1-\mathit{y}_1\,\mathit{y}_2)\,(\mathit{y}_1-\mathit{y}_2)) \,.
\end{align*}
We wish to compute the generating series
\[
 \sum _{\delta =0}^{\infty } 
\,{\Omega ^{(2)}_{x}}(\mathbf{T}(p^{\delta }))\cdot \Omega ^{(2)}_{y}(\mathbf{T}(p^{\delta }))\,X^{\delta }
\in\Q[x_0, x_1, x_2, y_0, y_1, y_2][\![X]\!].
\]
The answer  is given by the following multiplicative analogue  of Rankin's lemma in genus two:

\begin{theo}
 The following equality holds
\begin{align}
\label{Sxxyy}
 & \sum _{\delta =0}^{\infty } \,{\Omega ^{(2)}_{x}}(\mathbf{T}(p^{\delta }))\cdot
 {\Omega ^{(2)}_{y}}(\mathbf{T}(p^{\delta }))\,X^{\delta } = 
\\ & \nonumber 
-{\displaystyle \frac 
 {(p\,\mathit{x}_1 - \mathit{x}_2)\,(1-p\,\mathit{y}_1\,\mathit{y}_2)\,
   \mathit{x}_1\,\mathit{y}_1\,\mathit{y}_2}
 {p^{2}\,(1-\mathit{x}_1)\,(1-\mathit{x}_2)\,(\mathit{x}_1-\mathit{x}_2)\,
   (1-\mathit{y}_1)\,(1-\mathit{y}_2)\,(1-\mathit{y}_1\,\mathit{y}_2)\,
   (1-\mathit{x}_0\,\mathit{x}_1\,\mathit{y}_0\,\mathit{y}_1\,\mathit{y}_2\,X)}}
\\ & \nonumber 
+{\displaystyle \frac 
{\mathit{x}_2\,\mathit{y}_1\,(\mathit{x}_1-p\,\mathit{x}_2)\,(p\,\mathit{y}_1-\mathit{y}_2)}
{p^{2}\,(1-\mathit{x}_1)\,(1-\mathit{x}_2)\,(\mathit{x}_1-\mathit{x}_2)\,
  (1-\mathit{y}_1)\,(1-\mathit{y}_2)\,(\mathit{y}_1-\mathit{y}_2)\,
  (1-\mathit{x}_0\,\mathit{x}_2\,\mathit{y}_0\,\mathit{y}_1\,X)}} 
\\ & \nonumber 
+{\displaystyle \frac 
{\mathit{x}_2\,\mathit{y}_2\,(\mathit{x}_1-p\,\mathit{x}_2)(\mathit{y}_1-p\,\mathit{y}_2)}
{p^{2}\,(1-\mathit{x}_1)\,(1-\mathit{x}_2)\,(\mathit{x}_1-\mathit{x}_1)\,
  (1-\mathit{y}_1)\,(1-\mathit{y}_2)\,(\mathit{y1}-\mathit{y}_2)\,
  (1-\mathit{x}_0\,\mathit{y}_0\,\mathit{x}_2\,\mathit{y}_2\,X)}} 
\\ & \nonumber 
-{\displaystyle \frac 
{\mathit{x}_2\,\mathit{y}_1\,\mathit{y}_2\,(\mathit{x}_1-p\,\mathit{x}_2)\,
  (1-p\,\mathit{y}_1\,\mathit{y}_2)}
{p^{2}\,(1-\mathit{x}_1)\,(1-\mathit{x}_2)\,(\mathit{x}_1 -\mathit{x}_2)\,
  (1-\mathit{y}_1)\,(1-\mathit{y}_2)\,(1-\mathit{y}_1\,\mathit{y}_2)\,
(1-\mathit{x}_0\,\mathit{x}_2\,\mathit{y}_0\,\mathit{y}_1\,\mathit{y}_2\,X)}}
  \\ & \nonumber 
-{\displaystyle \frac 
{\mathit{x}_1\,(p\,\mathit{x}_1-\mathit{x}_2)\,(p-\mathit{y}_1\,\mathit{y}_2)}
{p^{2}\,(1-\mathit{x}_1)\,(1-\mathit{x}_2)\,(\mathit{x}_1-\mathit{x}_2)\,
  (1-\mathit{y}_1)\,(1-\mathit{y}_2)\,(1-\mathit{y}_1\,\mathit{y}_2)\,
  (1-\mathit{x}_0\,\mathit{x}_1\,\mathit{y}_0\,X)}}
  \\ & \nonumber 
-{\displaystyle \frac 
{\mathit{x}_1\,\mathit{x}_2\,\mathit{y}_1(1-p\,\mathit{x}_1\,\mathit{x}_2)\,
  (p\,\mathit{y}_1-\mathit{y}_2)}
{p^{2}\,(1-\mathit{x}_1)\,(1-\mathit{x}_2)\,(1-\mathit{x}_1\,\mathit{x}_2)\,
  (1-\mathit{y}_1)\,(1-\mathit{y}_2)\,(\mathit{y}_1-\mathit{y}_2)\,
  (1-\mathit{x}_0\,\mathit{x}_1\,\mathit{x}_2\,\mathit{y}_0\,\mathit{y}_1\,X)}}
\\ & \nonumber 
-{\displaystyle \frac 
{\mathit{x}_1\,\mathit{x}_2\,\mathit{y}_2\,(1-p\,\mathit{x}_1\,\mathit{x}_2)\,
  (\mathit{y}_1-p\,\mathit{y}_2)}
{p^{2}\,(1-\mathit{x}_1)\,(1-\mathit{x}_2)\,(1-\mathit{x}_1\,\mathit{x}_2)\,
  (1-\mathit{y}_1)\,(1-\mathit{y}_2)\,(\mathit{y}_1-\mathit{y}_2)\,
  (1-\mathit{x}_0\,\mathit{x}_1\,\mathit{x}_2\,\mathit{y}_0\,\mathit{y}_2\,X)}}
\\ & \nonumber 
+{\displaystyle \frac 
{\mathit{y}_1\,\mathit{y}_2\,(p-\mathit{x}_1\,\mathit{x}_2)\,(1-p\,\mathit{y}_1\,\mathit{y}_2)}
{p^{2}\,(1-\mathit{x}_1)\,(1-\mathit{x}_2)\,(1-\mathit{x}_1\,\mathit{x}_2)\,
  (1-\mathit{y}_1)\,(1-\mathit{y}_2)\,(1-\mathit{y}_1\,\mathit{y}_2)\,
  (1-\mathit{x}_0\,\mathit{y}_0\,\mathit{y}_1\,\mathit{y}_2\,X)}}
\\ & \nonumber 
+{\displaystyle \frac 
{\mathit{x}_1\,\mathit{x}_2\,(1-p\,\mathit{x}_1\,\mathit{x}_2)\,
  (p-\mathit{y}_1\,\mathit{y}_2)}
{p^{2}\,(1-\mathit{x}_1)\,(1-\mathit{x}_2)\,(1-\mathit{x}_1\,\mathit{x}_2)\,
  (1-\mathit{y}_1)\,(1-\mathit{y}_2)\,(1-\mathit{y}_1\,\mathit{y}_2)\,
  (1-\mathit{x}_0\,\mathit{x}_1\,\mathit{x}_2\,\mathit{y}_0\,X)}}
\end{align}
\end{theo}
\begin{align*}&
-{\displaystyle \frac 
{\mathit{x}_1\,\mathit{y}_1\,(p\,\mathit{x}_1-\mathit{x}_2)\,(p\,\mathit{y}_1-\mathit{y}_2)}
{p^{2}\,(1-\mathit{x}_1)\,(1-\mathit{x}_2)\,(\mathit{x}_1-\mathit{x}_2)\,
  (1-\mathit{y}_1)\,(1-\mathit{y}_2)\,(\mathit{y}_1-\mathit{y}_2)\,
  (1-\mathit{x}_0\,\mathit{x}_1\,\mathit{y}_0\,\mathit{y}_1\,X)}}
\\ & \nonumber 
+{\displaystyle \frac 
{\mathit{x}_1\,\mathit{y}_2\,(p\,\mathit{x}_1-\mathit{x}_2)\,(\mathit{y}_1-p\,\mathit{y}_2)}
{p^{2}\,(1-\mathit{x}_1)\,(1-\mathit{x}_2)\,(\mathit{x}_1-\mathit{x}_2)\,
  (1-\mathit{y}_1)\,(1-\mathit{y}_2)\,(\mathit{y}_1-\mathit{y}_2)\,
(1-\mathit{x}_0\,\mathit{x}_1\,\mathit{y}_0\,\mathit{y}_2\,X)}}
\\ & \nonumber 
-{\displaystyle \frac 
{\mathit{x}_2\,(\mathit{x}_1-p\,\mathit{x}_2)\,(p-\mathit{y}_1\,\mathit{y}_2)}
{p^{2}\,(1-\mathit{x}_1)\,(1-\mathit{x}_2)\,(\mathit{x}_1-\mathit{x}_2)\,
  (1-\mathit{y}_1)\,(1-\mathit{y}_2)\,(1-\mathit{y}_1\,\mathit{y}_2)\,
  (1-\mathit{x}_0\,\mathit{x}_2\,\mathit{y}_0\,X)}}
\\ & \nonumber 
+{\displaystyle \frac 
{\mathit{x}_1\,\mathit{x}_2\,\mathit{y}_1\,\mathit{y}_2\,(1-p\,\mathit{x}_1\,\mathit{x}_2)\,
  (1-p\,\mathit{y}_1\,\mathit{y}_2)}
{p^{2}\,(1-\mathit{x}_1)\,(1-\mathit{x}_2)\,(1-\mathit{x}_1\,\mathit{x}_2)\,
  (1-\mathit{y}_1)\,(1-\mathit{y}_2)\,(1-\mathit{y}_1\,\mathit{y}_2)\,
(1-\mathit{x}_0\,\mathit{x}_1\,\mathit{x}_2\,\mathit{y}_0\,\mathit{y}_1\,\mathit{y}_2\,X)}}
\\ & \nonumber 
+{\displaystyle \frac 
{(p-\mathit{x}_1\,\mathit{x}_2)\,(p-\mathit{y}_1\,\mathit{y}_2)}
{p^{2}\,(1-\mathit{x}_1)\,(1-\mathit{x}_2)\,(1-\mathit{x}_1\,\mathit{x}_2)\,
  (1-\mathit{y}_1)\,(1-\mathit{y}_2)\,(1-\mathit{y}_1\,\mathit{y}_2)\,
  (1-\mathit{x}_0\,\mathit{y}_0\,X)}}
\\ & \nonumber 
-{\displaystyle \frac 
{\mathit{y}_1\,(p-\mathit{x}_1\,\mathit{x}_2)\,(p\,\mathit{y}_1 - \mathit{y}_2)}
{p^{2}\,(1-\mathit{x}_1)\,(1-\mathit{x}_2)\,(1-\mathit{x}_1\,\mathit{x}_2)\,
  (1-\mathit{y}_1)\,(1-\mathit{y}_2)\,(\mathit{y}_1-\mathit{y}_2)\,
  (1-\mathit{x}_0\,\mathit{y}_0\,\mathit{y}_1\,X)}}
\\ & \nonumber 
-{\displaystyle \frac 
{\mathit{y}_2\,(p-\mathit{x}_1\,\mathit{x}_2)\,(\mathit{y}_1-p\,\mathit{y}_2)}
{p^{2}\,(1-\mathit{x}_1)\,(1-\mathit{x}_2)\,(1-\mathit{x}_1\,\mathit{x}_2)\,
  (1-\mathit{y}_1)\,(1-\mathit{y}_2)\,(\mathit{y}_1-\mathit{y}_2)\,
  (1-\mathit{x}_0\,\mathit{y}_0\,\mathit{y}_2\,X)}}
\end{align*}

\begin{rema}[On the denominator of (\ref{Sxxyy})]
One finds using a computer  that the polynomials not depending on  $X$
in the denominators of (\ref{Sxxyy}) cancel after the simplification 
in the ring $\Q[x_0, x_1, x_2, y_0, y_1, y_2][\![X]\!]$, 
so that the common denominator becomes
\begin{align*} &
(1 - \mathit{x}_0\,\mathit{y}_0 X)\,
(1 - \mathit{x}_0\,\mathit{y}_0\,\mathit{x}_1\,X)\,
(1 - \mathit{x}_0\,\mathit{y}_0\,\mathit{y}_1\,X)\,
(1 - \mathit{x}_0\,\mathit{y}_0\,\mathit{x}_2\,X)\,
(1 - \mathit{x}_0\,\mathit{y}_0\,\mathit{y}_2\,X)\,
 \\ &
(1 - \mathit{x}_0\,\mathit{y}_0\,\mathit{x}_1\,\mathit{y}_1\,X)\,
(1 - \mathit{x}_0\,\mathit{y}_0\,\mathit{x}_1\,\mathit{x}_2\,X)\,
(1 - \mathit{x}_0\,\mathit{y}_0\,\mathit{x}_1\,\mathit{y}_2\,X)\,
(1 - \mathit{x}_0\,\mathit{y}_0\,\mathit{y}_1\,\mathit{x}_2\,X)\,
 \\ &
(1 - \mathit{x}_0\,\mathit{y}_0\,\mathit{y}_1\,\mathit{y}_2\,X)\,
(1 - \mathit{x}_0\,\mathit{y}_0\,\mathit{x}_2\,\mathit{y}_2\,X)\,
(1 - \mathit{x}_0\,\mathit{y}_0\,\mathit{x}_1\,\mathit{y}_1\,\mathit{x}_2\,X)\,
(1 - \mathit{x}_0\,\mathit{y}_0\,\mathit{x}_1\,\mathit{y}_1\,\mathit{y}_2\,X)\,
 \\ &
(1 - \mathit{x}_0\,\mathit{y}_0\,\mathit{x}_1\,\mathit{x}_2\,\mathit{y}_2\,X)\,
(1 - \mathit{x}_0\,\mathit{y}_0\,\mathit{y}_1\,\mathit{x}_2\,\mathit{y}_2\,X)\,
(1 - \mathit{x}_0\,\mathit{y}_0\,\mathit{x}_1\,\mathit{y}_1\,\mathit{x}_2\,\mathit{y}_2\, X).
\end{align*}
 
\end{rema}
\begin{rema}[Comparison with $g=1$]\label{Rn1} 
It turns out by direct computation that the numerator of the rational fraction
(\ref{Sxxyy})
is a product of the factor $1-\mathit{x}_0^{2}\mathit{y}_0^{2}\,
\mathit{x}_1\mathit{y}_1\mathit{x}_2\mathit{y}_2X^{2}$
by a polynomial of degree $12$ in $X$ with coefficients in $\Q[x_0, x_1, x_2, y_0, y_1, y_2]$
with the constant term equal to 1 and the leading term  
$$
\frac{x_0^{12}y_0^{12}x_1^6x_2^6y_1^6y_2^6}{p^2}X^{12}.
$$
Moreover, the factor of degree 12 does not contain terms of degree 1 and 11 in $X$.
The factor of  degree 2 in $X$ is very  similar to one in the case $g=1$ (this series was studied and used in \cite{Ma-Pa77}
): 
\begin{align*}
 & \sum _{\delta =0}^{\infty }\,
{\Omega ^{(1)}_{x}}(\mathbf{T}(p^{\delta }))\cdot 
{\Omega ^{(1)}_{y}}(\mathbf{T}(p^{\delta }))\,X^{\delta } = 
\sum _{\delta =0}^{\infty }
{\displaystyle \frac 
{\mathit{x}_0^{\delta }\, (1-\mathit{x}_1^{(1+\delta )})}{1-\mathit{x}_1}} 
\cdot
{\displaystyle \frac 
{\mathit{y}_0^{\delta }\, (1-\mathit{y}_1^{(1+\delta )})}{1-\mathit{y}_1}} 
X^{\delta }
\\ &  =
{\displaystyle \frac {1}
{(1-\mathit{x}_1)\,(1-\mathit{y}_1)\,(1-\mathit{x}_0\,\mathit{y}_0\,X)}}
 - {\displaystyle \frac {\mathit{y}_1}
{(1-\mathit{x}_1)\,(1-\mathit{y}_1)\,(1-\mathit{x}_0\,\mathit{y}_0\,\mathit{y}_1\,X)}} 
 \\ &
- {\displaystyle \frac {\mathit{x}_1}
{(1-\mathit{x}_1)\,(1-\mathit{y}_1)\,(1-\mathit{x}_0\,\mathit{y}_0\,\mathit{x}_1\,X)}}
+ {\displaystyle \frac {\mathit{x}_1\,\mathit{y}_1}
{(1-\mathit{x}_1)\,(1-\mathit{y}_1)\,(1-\mathit{x}_0\,\mathit{y}_0\,\mathit{x}_1\,\mathit{y}_1\,X)}}
\\ & =
{\displaystyle \frac {1-\mathit{x}_0^{2}\,\mathit{y}_0^{2}\,\mathit{x}_1\,\mathit{y}_1\,X^{2}}
{(1-\mathit{x}_0\,\mathit{y}_0\,\mathit{x}_1\,\mathit{y}_1\,X)\,
 (1-\mathit{x}_0\,\mathit{y}_0\,\mathit{x}_1\,X)\,
 (1-\mathit{x}_0\,\mathit{y}_0\,\mathit{y}_1\,X)\,
 (1-\mathit{x}_0\,\mathit{y}_0\,X)}}\,.
\end{align*}
\end{rema}
\subsection{Symmetric square generating series in genus 2} 
Using the same method, one can evaluate the symmetric square generating series
and the cubic generating series of higher genus.
Note that this series, written here in spherical variables $x_0, x_1, x_2$
is different from one studied by Andrianov-Kalinin, and has the form: 

\begin{align*} & 
{\displaystyle \sum _{\delta =0}^{\infty }} \,
{\Omega^{(2)} _{x}}(\mathbf{T}(p^{2\delta }))X^{\delta }=
\\ &
\frac{(1
+  \mathit{x}_0^{2}\,\mathit{x}_1\,X 
+  \mathit{x}_0^{2}\,\mathit{x}_2\,X
+ 2\mathit{x}_0^{2}\,\mathit{x}_1\,\mathit{x}_2\,X
+  \mathit{x}_0^{2}\,\mathit{x}_1\,\mathit{x}_2^{2}\,X
+  \mathit{x}_0^{2}\,\mathit{x}_1^{2}\,\mathit{x}_2\,X 
+  \mathit{x}_0^{4}\,\mathit{x}_1^{2}\,\mathit{x}_2^{2}\,X^{2})
(1-\frac{\mathit{x}_0^{2}\mathit{x}_1\mathit{x}_2}{p} X)}
{(1 - \mathit{x}_0^{2}\,\mathit{x}_1^{2}\,\mathit{x}_2^{2}\,X) 
 (1 - \mathit{x}_0^{2}\,\mathit{x}_1^{2}\,X)
 (1 - \mathit{x}_0^{2}\,\mathit{x}_2^{2}\,X)
 (1 - \mathit{x}_0^{2}\,X)}\,.
\end{align*}

\subsection{Cubic generating series in genus 2}
The cubic generating series of higher genus, written here in spherical variables $x_0, x_1, x_2$
 has the form: 

\begin{align*} &
{\displaystyle \sum _{\delta =0}^{\infty }} \,{\Omega^{(2)} _{x}}(
\mathbf{T}(p^{3\,\delta }))\,X^{\delta }= 
p^{-1} ( - p + \mathit{x}_0
^{6}\,\mathit{x}_1^{4}\,\mathit{x}_2^{2}\,X^{2} + \mathit{x}_0^{6}\,
\mathit{x}_1^{2}\,\mathit{x}_2^{4}\,X^{2} + 2\,\mathit{x}_0^{6}\,
\mathit{x}_1^{2}\,\mathit{x}_2^{3}\,X^{2} \\ &
\mbox{} - p\,\mathit{x}_0^{6}\,\mathit{x}_1^{4}\,\mathit{x}_2^{4}\,X
^{2} - p\,\mathit{x}_0^{6}\,\mathit{x}_1^{2}\,\mathit{x}_2^{4}\,X^{2
} - 2\,p\,\mathit{x}_0^{3}\,\mathit{x}_1^{2}\,\mathit{x}_2\,X + 
\mathit{x}_0^{6}\,\mathit{x}_1\,\mathit{x}_2^{3}\,X^{2} \\ &
\mbox{} + \mathit{x}_0^{6}\,\mathit{x}_1^{3}\,\mathit{x}_2\,X^{2} + 
\mathit{x}_0^{6}\,\mathit{x}_1^{3}\,\mathit{x}_2^{5}\,X^{2} + 
\mathit{x}_0^{6}\,\mathit{x}_1^{5}\,\mathit{x}_2^{3}\,X^{2} + 3\,
\mathit{x}_0^{6}\,\mathit{x}_1^{3}\,\mathit{x}_2^{3}\,X^{2} \\ &
\mbox{} + \mathit{x}_0^{6}\,\mathit{x}_1^{2}\,\mathit{x}_2^{2}\,X^{2
} + 2\,\mathit{x}_0^{6}\,\mathit{x}_1^{3}\,\mathit{x}_2^{2}\,X^{2}
 - p\,\mathit{x}_0^{3}\,\mathit{x}_1^{2}\,X - p\,\mathit{x}_0^{3}\,
\mathit{x}_2^{2}\,X - p\,\mathit{x}_0^{6}\,\mathit{x}_1^{4}\,
\mathit{x}_2^{2}\,X^{2} \\ &
\mbox{} - 2\,p\,\mathit{x}_0^{3}\,\mathit{x}_1\,\mathit{x}_2^{2}\,X
 - p\,\mathit{x}_0^{6}\,\mathit{x}_1^{2}\,\mathit{x}_2^{2}\,X^{2} + 
\mathit{x}_0^{3}\,\mathit{x}_1^{2}\,\mathit{x}_2^{2}\,X + \mathit{x}_0^{3}\,\mathit{x}_1\,\mathit{x}_2\,X \\ &
\mbox{} - p\,\mathit{x}_0^{6}\,\mathit{x}_1^{2}\,\mathit{x}_2^{3}\,X
^{2} - p\,\mathit{x}_0^{6}\,\mathit{x}_1^{3}\,\mathit{x}_2^{2}\,X^{2
} - p\,\mathit{x}_0^{3}\,\mathit{x}_1^{2}\,\mathit{x}_2^{3}\,X - p\,
\mathit{x}_0^{3}\,\mathit{x}_1^{3}\,\mathit{x}_2^{2}\,X \\ &
\mbox{} - 2\,p\,\mathit{x}_0^{3}\,\mathit{x}_1^{2}\,\mathit{x}_2^{2}
\,X - p\,\mathit{x}_0^{3}\,\mathit{x}_1^{3}\,\mathit{x}_2\,X + 
\mathit{x}_0^{3}\,\mathit{x}_1^{2}\,\mathit{x}_2\,X + \mathit{x}_0^{9
}\,\mathit{x}_1^{4}\,\mathit{x}_2^{4}\,X^{3} \\ &
\mbox{} - 2\,p\,\mathit{x}_0^{6}\,\mathit{x}_1^{3}\,\mathit{x}_2^{3}
\,X^{2} - 2\,p\,\mathit{x}_0^{3}\,\mathit{x}_1\,\mathit{x}_2\,X + 
\mathit{x}_0^{9}\,\mathit{x}_1^{4}\,\mathit{x}_2^{5}\,X^{3} + 
\mathit{x}_0^{9}\,\mathit{x}_1^{5}\,\mathit{x}_2^{4}\,X^{3} \\ &
\mbox{} - p\,\mathit{x}_0^{6}\,\mathit{x}_1^{3}\,\mathit{x}_2^{4}\,X
^{2} + \mathit{x}_0^{3}\,\mathit{x}_1\,\mathit{x}_2^{2}\,X + 
\mathit{x}_0^{9}\,\mathit{x}_1^{5}\,\mathit{x}_2^{5}\,X^{3} + 
\mathit{x}_0^{6}\,\mathit{x}_1^{4}\,\mathit{x}_2^{4}\,X^{2} \\ &
\mbox{} - p\,\mathit{x}_0^{6}\,\mathit{x}_1^{4}\,\mathit{x}_2^{3}\,X
^{2} - p\,\mathit{x}_0^{3}\,\mathit{x}_1\,\mathit{x}_2^{3}\,X - p\,
\mathit{x}_0^{3}\,\mathit{x}_2\,X - p\,\mathit{x}_0^{3}\,\mathit{x}_1
\,X + 2\,\mathit{x}_0^{6}\,\mathit{x}_1^{4}\,\mathit{x}_2^{3}\,X^{2}
 \\ &
\mbox{} + 2\,\mathit{x}_0^{6}\,\mathit{x}_1^{3}\,\mathit{x}_2^{4}\,X
^{2})/(( 1 - \mathit{x}_0^{3}\,X)( 1 - \mathit{x}_0^{3}
\mathit{x}_1^{3}\,X)( 1 - \mathit{x}_0^{3}\mathit{x}_2^{3}\,X)
(  1 - \mathit{x}_0^{3}\mathit{x}_1^{3}\mathit{x}_2^{3}X))\,.
\end{align*}

\section{Proofs:  formulas for the  Hecke operators of $Sp_g$}
\subsection{Satake's spherical map $\Omega$}
Our result is based on the use of the Satake spherical map $\Omega$, by
applying the spherical map $\Omega$ 
to elements $\mathbf{T}(p^\delta )  \in {\mathcal L}_\Z$ of Hecke ring
${\mathcal L}_\Z=\Z[ \mathbf{T}(p), \mathbf{T}_1(p^2), \cdots, \mathbf{T}_{n}(p^2)]$
 for the symplectic 
group, see  \cite{An87} chapter 3.  

\begin{itemize}
\item{Case $\mathbf{T}_1(p^2)$.} 
In genus 2 (in spherical variables $x_0, x_1, x_2$), we obtain 
using  Andrianov's formulas:
\begin{align*}&
\Omega(\mathbf{T}_1(p^2)) = {\displaystyle \frac {\mathit{x}_0^{2}\,((
\mathit{x}_1^{2}\,\mathit{x}_2 + \mathit{x}_1\,\mathit{x}_2^{2})\,p^{
2} + \mathit{x}_1\,\mathit{x}_2\,p^{2} - \mathit{x}_1\,\mathit{x}_2
 + (\mathit{x}_1 + \mathit{x}_2)\,p^{2})}{p^{3}}}\,.
\end{align*}
\item{Cases $\mathbf{T}_2(p^2)=[\mathbf{p}]_2$ and $\mathbf{T}(p)$}
$$
\Omega(\mathbf{T}_2(p^2)) =  
 \frac 
{\mathit{x}_0^{2}\,
\mathit{x}_1\,\mathit{x}_2}{p^{3}} ,\ \ 
\Omega(\mathbf{T}(p)) = \mathit{x}_0(1+\mathit{x}_1)(1+\mathit{x}_2)\,.
$$
\end{itemize}


\subsection{Use of Andrianov's generating series in genus 2}
We refer to \cite{An87}, p.164, (3.3.75)
for the following celebrated summation formula:
\begin{align}
  \sum _{\delta =0}^{\infty } \,{\Omega }^{(2)}(\mathbf{T}(p^{\delta }))\,X^{\delta }=
{\displaystyle \frac 
{1-\frac{\mathit{x}_0^{2}\,\mathit{x}_1\,\mathit{x}_2}{p}\,X^{2}}
{
 (1-\mathit{x}_0\,X)\,
 (1-\mathit{x}_0\,\mathit{x}_1\,X)
 (1-\mathit{x}_0\,\mathit{x}_2\,X)\,
 (1-\mathit{x}_0 \,\mathit{x}_1\,\mathit{x}_2\,X)\,
}}\,,
\end{align}
which gives after development and simplification the following formula 
\begin{align*} 
{\Omega }^{(2)}(\mathbf{T}&(p^{\delta })) = 
  p^{-1}\,\mathit{x}_0^{ \delta }
( p\,\mathit{x}_1^{(3+\delta)}\,\mathit{x}_2 
- p\,\mathit{x}_1^{(2+\delta)} 
- p\,\mathit{x}_1^{(3+\delta)}\,\mathit{x}_2^{(2+\delta)} 
+ p\,\mathit{x}_1^{(2+\delta)}\,\mathit{x}_2^{(3+\delta)} 
\\ & \quad 
- p\,\mathit{x}_1\,\mathit{x}_2^{(3+\delta)} 
+ p\,\mathit{x}_2^{(2+\delta)} 
+ p\,\mathit{x}_1 
- p\,\mathit{x}_2 
- \mathit{x}_1^{(2+\delta)}\,\mathit{x}_2^{2} 
+ \mathit{x}_1^{(1+\delta)}\,\mathit{x}_2 
\\ & \quad 
+ \mathit{x}_1^{(2+\delta)}\,\mathit{x}_2^{(1+\delta)} 
- \mathit{x}_1^{(1+\delta)}\,\mathit{x}_2^{(2+\delta)} 
+ \mathit{x}_1^{2}\,\mathit{x}_2^{(2+\delta)} 
- \mathit{x}_1\,\mathit{x}_2^{(1+\delta)} 
- \mathit{x}_1^{2}\,\mathit{x}_2 
+ \mathit{x}_1\,\mathit{x}_2^{2})/ 
\\ & \quad 
((1-\mathit{x}_1)\,(1-\mathit{x}_2)\,(1-\mathit{x}_1\,\mathit{x}_2)\,(\mathit{x}_1-\mathit{x}_2))\,.
\end{align*}

Then we use  two groups of variables: $x_0,\dots ,x_n$ and $y_0,\dots ,y_n$ 
in  two copies  ${\Omega _{x}}$, ${\Omega _{y}}$ of the sperical map, 
in order to treat the tensor product of two local Hecke algebras. 

Next, in order to carry out the summation of the series
$$
\sum _{\delta =0}^{\infty } 
\,{\Omega ^{(2)}_{x}}(\mathbf{T}(p^{\delta }))\cdot \Omega ^{(2)}_{y}(\mathbf{T}(p^{\delta }))\,X^{\delta }
$$
on a computer, we used a subdivision of each summand (over $\delta$) into smaller parts.
These parts correspond to symbolic monomials in $x_1^\delta, y_1^\delta, 
x_2^\delta, y_2^\delta, 
(x_1x_2)^\delta,  (y_1y_2)^\delta$.



\def\Ta
{\mathbf{T}(p)}
\def\Tb
{\mathbf{T}_1(p^2)}
\def\Tc
{[\mathbf{p}] }
\def\Txa
{\mathbf{T}(p)\otimes 1}
\def\Txb
{\mathbf{T}_1(p^2)\otimes 1}
\def\Txc
{[\mathbf{p}] \otimes 1}
\def\Tya
{1\otimes \mathbf{T}(p)}
\def\Tyb
{1\otimes \mathbf{T}_1(p^2)}
\def\Tyc
{1 \otimes [\mathbf{p}]}
\def\Taa
{\mathbf{T}(p)^2}
\def\Txaa
{\mathbf{T}(p)^2\otimes 1}
\def\Tba
{\mathbf{T}_1(p^2)^2}
\def\Txba
{\mathbf{T}_1(p^2)^2\otimes 1}
\def\Tca
{[\mathbf{p}]^2}
\def\Txca
{[\mathbf{p}]^2\otimes 1}
\def\Tyaa
{1\otimes \mathbf{T}(p)^2}
\def\Tyba
{1\otimes \mathbf{T}_1(p^2)^2}
\def\Tyca
{1 \otimes [\mathbf{p}]^2}
\def\Tab
{\mathbf{T}(p)^3} 
\def\Txab
{\mathbf{T}(p)^3\otimes 1}
\def\Tbb
{\mathbf{T}_1(p^2)^3}
\def\Txbb
{\mathbf{T}_1(p^2)^3\otimes 1}
\def\Tcb
{ [\mathbf{p}]^3}
\def\Txcb
{1 \otimes [\mathbf{p}]^3}
\def\Tyab
{1\otimes \mathbf{T}(p)^3}
\def\Tybb
{1\otimes \mathbf{T}_1(p^2)^3}
\def\Tycb
{1 \otimes [\mathbf{p}]^3}
\subsection{Rankin's Lemma of genus  2 (compare with [Jia96])}

Let us compute the series
\[
D^{(1,1)}_p(X)= \sum _{\delta =0}^{\infty } 
\mathbf{T}(p^{\delta })\otimes\mathbf{T}(p^{\delta })\,X^{\delta }
\in{\mathcal L}_{2,\Z}\otimes{\mathcal L}_{2,\Z}[\![X]\!]
\]
in terms of the generators of Hecke's algebra
${\mathcal L}_{2,\Z}\otimes{\mathcal L}_{2,\Z}$
given by the following operators:
\begin{align*}&
\mathbf{T}(p)\otimes 1,
\mathbf{T}_1(p^2)\otimes 1,
[\mathbf{p}]\otimes 1, 
1\otimes \mathbf{T}(p),
1\otimes \mathbf{T}_1(p^2),
1 \otimes [\mathbf{p}]\in{\mathcal L}_{2,\Z}
\otimes{\mathcal L}_{2,\Z}[\![X]\!].
\end{align*}

\begin{theo}[October 2006] 
\label{ThRS} 
For $g=2$, we have the following explicit representation 

\begin{align*}&
D^{(1,1)}_p(X)=\sum _{\delta =0}^{\infty } 
\mathbf{T}(p^{\delta })\otimes\mathbf{T}(p^{\delta })\,X^{\delta }
=\frac{(1-p^6 [\mathbf{p}] \otimes [\mathbf{p}] X^2)\cdot R(X)}{S(X)}, 
\mbox{ where} 
\\ & 
R(X), S(X)
%
\in{\mathcal L}_{2,\Z}\otimes
{\mathcal L}_{2,\Z}[X] 
\end{align*}
are given by the  equalities  (\ref{EqR}) and (\ref{EqS}):
\begin{align}\label{EqR} 
R(X)&=1+r_2X^2+\cdots +r_{10}X^{10}+r_{12}X^{12} \in{\mathcal L}_{2,\Z}\otimes {\mathcal L}_{2,\Z}[X] \mbox{ with } r_{1}=r_{11}=0,
\\ 
\label{EqS} 
S(X)&=1+s_1X+\cdots +s_{16}X^{16}
\\ &\nonumber
= 1 - (\Ta \otimes \Ta) X+\cdots +(p^6\Tc\otimes \Tc)^{8}X^{16}\in{\mathcal L}_{2,\Z}\otimes{\mathcal L}_{2,\Z}[X],
\end{align}
 with $r_i$ and $s_i$  given in Appendix. 
Moreover, there is an easy  functional equation   (similar to \cite{An87}, p.164, (3.3.79)):   
$$
s_{16-i}=(p^6\Tc\otimes \Tc)^{8-i}s_i \ \ \ (i=0, \cdots, 8).
$$ 
 \end{theo}

\begin{rema}[Comparison with the case  $g=1$ (in terms of Hecke operators)]
The corresponding result in the case  $g=1$ written in terms of Hecke operators, looks as follows (see \ref{Rn1}): 
\hskip-.5cm\begin{align*}
& \sum _{\delta =0}^{\infty }\,
\mathbf{T}(p^{\delta })\otimes 
\mathbf{T}(p^{\delta })\,X^{\delta } = 
(1-
p^2\Tc\otimes\Tc
X^{2})/
 \\ &
(1 - 
 \Ta\otimes\Ta
X 
+ \big(p
 (\Ta^2\otimes\Tc+ \Tc\otimes\Ta^2)
- 2p^2
\Tc\otimes\Tc
\big)X^{2} 
 \\ &
 - p^{2}
\Ta\Tc\otimes\Ta\Tc
X^{3}
 + p^{4}
\Tc^2\otimes\Tc^2 
X^{4}
).
\end{align*}
Indeed this follows directly from Remark \ref{Rn1}
\end{rema}


\section{Applications to $L$-functions and motives for $\Sp_n$
}
\subsubsection*{The Fourier expansion of a Siegel modular form.
}

Let  $\ds f=\sum_{{\mathcal T}\in B_n} a({\mathcal T})q^{\mathcal T}\in {\mathcal M}_k^n$ be a Siegel modular form of weight $k$ and of genus $n$
on the Siegel upper-half plane $\HH_n=\{z\in {\mathrm M}_n(\CC) \ |\ \Im(z) >0 \}$.

\

The formal  {\it  Fourier expansion }  of $f$ uses the symbol
$$
q^{\mathcal T}=\exp(2\pi i {\rm tr}({\mathcal T}z)) = 
 \prod_{i=1}^nq_{ii}^{{\mathcal T}_{ii}}\prod_{i<j}q_{ij}^{2{\mathcal T}_{ij}}
$$
$ \in \dbC[\![q_{11}, \dots , q_{nn}]\!][q_{ij}, \ q_{ij}^{-1}]_{i,j=1, \cdots, m}$, where
$q_{ij}= \exp(2\pi (\sqrt {-1} z_{i,j}))$, 
and
${\mathcal T}$ is in the   semi-group
 $
B_n =\{{\mathcal T}={}^t{\mathcal T}\ge 0\vert {\mathcal T}\hbox{ half-integral}\}. 
$
\subsubsection*{Satake parameters of an eigenfunction of Hecke operators}
Suppose that $f\in\mathcal{M}_k^n$ is an eigenfunction of all  {\it Hecke operators}
$f \longmapsto f|T$, $T\in {\mathcal L}_{n,p}$ for all primes $p$, 
hence $f|T = \lambda_f(T) f$. 

Then all the numbers 
$\lambda_f(T) \in \CC$ define a homomorphism
$\lambda_f : {\mathcal L}_{n,p} \longrightarrow \CC$ 
given by a  $(n + 1)$-tuple of complex numbers
$
    (\al_0, \al_1, \cdots, \al_n) =
     (\CC^{\times})^{n+1}  
$
(the Satake parameters of $f$).

One has $
     \al_0^2 \al_1 \cdots \al_n=  p^{kn-n(n+1)/2}.  $

For another Siegel modular form, eigenfunction of Hecke operators
  $g \in{\mathcal M}_l^n$ consider the corresponding homomorphism
$\lambda_g : {\mathcal L}_{n,p} \longrightarrow \CC$ 
 given by its Satake parameters  $ (\beta_0, \beta_1, \cdots, \beta_n)$ of $g$, and let
$\lambda_f\otimes\lambda_g : {\mathcal L}_{n,p}\otimes{\mathcal L}_{n,p} \longrightarrow \CC$.

\subsubsection*{  $L$-functions, functional equation and motives for   $\Sp_n$ (see \cite{Pa94}, \cite{Y})
}
One defines
\begin{itemize}
\item
$\ds
    Q_{f,p}(X)  =  (1 - \al_0 X) \prod_{r=1}^n
    \prod_{1 \leq i_1 < \cdots < i_r \leq n} (1 - \al_0 \al_{i_1}
    \cdots \al_{i_r} X),  $  

\item
 $\ds R_{f,p}(X)  =  (1 -  X)\prod_{i=1}^n (1 - \al_i^{-1} X)
    (1 - \al_i X) \in \QQ[ \al_0^{\pm 1}, \cdots, \al_n^{\pm 1}][X].$
\end{itemize}

Then the spinor  $L$-function $L(Sp(f), s, \chi)$ and the standard $L$-function $L(St(f), s, \chi)$ of $f$ (for $s\in \CC$, and for
all Dirichlet characters $\chi$) are defined as the Euler products  \\
\begin{itemize}
\item
 $ \ds L(Sp(f), s, \chi) = \prod_{p}Q_{f,p}(\chi(p)p^{-s})^{-1} 
$
\item
$\ds
 L(St(f), s, \chi)    =    \prod_{p}R_{f,p}(\chi(p)p^{-s})^{-1}
$
\end{itemize}
\subsubsection*{Motivic $L$-functions
}
Following  \cite{Pa94} and \cite{Y}, these functions are conjectured to be 
motivic for all $k>n$:
$$
   L(Sp(f), s, \chi) =L(M(Sp(f))(\chi), s),\, L(St(f), s, \chi) =L(M(St(f))(\chi), s),
$$
and the motives $M(Sp(f))$ and $M(St(f))$ are {\it pure} if $f$ is a genuine cusp form (not coming from a lifting of a smaller genus):
\begin{itemize}
{}\item
$M(Sp(f))$ is a motive over $\Q$ with coefficients in $\Q(\lambda_f(n))_{n\in\N}$
of  rank $2^n$, of weight $w=kn-n(n+1)/2$, and of  Hodge type
$\oplus_{p,q} H^{p,q}$, with 
\begin{align}\label{Hodge} &
p=(k-i_1)+(k-i_2)+\cdots+(k-i_r),\\ \nonumber &
q=(k-j_1)+(k-j_2)+\cdots+(k-i_s),\mbox{ where }r+s=n,\\ \nonumber &
1\le i_1<i_2<\cdots<i_r\le n, 1\le j_1<j_2<\cdots<j_s\le n,\\ \nonumber & 
\{ i_1,\cdots,i_r\}\cup \{ j_1,\cdots,i_s\}=\{ 1,2,\cdots,n\};
\end{align}
{}\item
$M(St(f))$ is a motive over $\Q$ with coefficients in $\Q(\lambda_f(n))_{n\in\N}$
of  rank $2n+1$, of weight $w=0$, and of  Hodge type
$H^{0,0}\oplus_{i=1}^n (H^{-k+i,k-i}\oplus H^{k-i, -k+i})$.
\end{itemize}
\subsubsection*{A functional equation 
}
Following  general Deligne's conjecture \cite{D} 
on the motivic $L$-functions,
the $L$-function satisfy a functional equation  
determined by the Hodge structure of a motive:
$$
\Lambda(Sp(f), kn-n(n+1)/2+1-s)=  \varepsilon(f) \Lambda(Sp(f), s),  \mbox{ where}
$$ 
$
\Lambda(Sp(f), s)=\Gamma_{n,k}(s)L(Sp(f), s),  \varepsilon(f)=(-1)^{k2^{n-2}},
$ \\ \  \\ {}
 $\Gamma_{1,k}(s)=\Gamma_\CC(s)=2(2\pi)^{-s}\Gamma(s)$, $\Gamma_{2,k}(s)=\Gamma_\CC(s)\Gamma_\CC(s-k+2)$, and\\
$\Gamma_{n,k}(s)=\prod_{p<q}\Gamma_\CC(s-p)\Gamma^{a_+}_\RR(s-(w/2))\Gamma_\RR(s+1-(w/2))^{a_-}$ for some non-negative integers
$a_+$ and $a_-$, with
$a_++a_-=w/2$, and $\Gamma_\RR(s)=\pi^{-s/2}\Gamma(s/2)$.

\subsubsection*{Motive of the  Rankin product of genus $g=2$}
Let $f$ and $g$ be two  Siegel cusp eigenforms of weights $k$ and $l$, $k>l$, and
let $M(Sp(f))$ and $M(Sp(g))$ be the spinor motives of  $f$ and $g$.
Then  $M(Sp(f))$ is a motive over $\Q$ with coefficients in $\Q(\lambda_f(n))_{n\in\N}$
of rank  $4$, of weight $w=2k-3$, and of  Hodge type
$H^{0,2k-3}\oplus H^{k-2,k-1}\oplus H^{k-1,k-2}\oplus H^{2k-3,0}$, 
and 
$M(Sp(g))$ is a motive over $\Q$ with coefficients in $\Q(\lambda_g(n))_{n\in\N}$
of rank $4$, of weight $w=2l-3$, and of  Hodge type
$H^{0,2l-3}\oplus H^{l-2,l-1}\oplus H^{l-1,l-2}\oplus H^{2l-3,0}$.

The tensor product $M(Sp(f))\otimes M(Sp(g))$ 
is a motive over $\Q$ with coefficients in  $\Q(\lambda_f(n), \lambda_g(n))_{n\in\N}$
of rank $16$, of weight $w=2k+2l-6$, and of  Hodge type

\begin{align*}&
H^{0,2k+2l-6}\oplus H^{l-2,2k+l-4}\oplus H^{l-1,2k+l-5}\oplus H^{2l-3,2k-3}
\\ &
H^{k-2,k+2l-4}\oplus H^{k+l-4,k+l-2}\oplus H_+^{k+l-3,k+l-3}\oplus H^{k+2l-5,k-1}
\\ &
H^{k-1,k+2l-5}\oplus H_-^{k+l-3,k+l-3}\oplus H^{k+l-2,k+l-4}\oplus H^{k+2l-4,k-2}
\\ &
H^{2k-3,2l-3}\oplus H^{2k+l-5,l-1}\oplus H^{2k+l-4,l-2}\oplus H^{2k+2l-6,0}.
\end{align*}

\subsubsection*{Motivic  $L$-functions: analytic properties
}
Following  Deligne's conjecture \cite{D} 
on motivic $L$-functions, applied 
for a Siegel cusp eigenform $F$ for the Siegel modular group $\Sp_4(\Z)$ of genus $n=4$ and  of weight  $k>5$, 
one has
$\Lambda(Sp(F), s)=\Lambda(Sp(F), 4k-9-s)$,  
where
\begin{align*}&
\Lambda(Sp(F), s)=\Gamma_\CC(s)\Gamma_\CC(s-k+4)\Gamma_\CC(s-k+3)\Gamma_\CC(s-k+2)\Gamma_\CC(s-k+1)
\\ & \times
\Gamma_\CC(s-2k+7)\Gamma_\CC(s-2k+6)\Gamma_\CC(s-2k+5)
L(Sp(F), s),
\end{align*}
(compare this functional equation with that given in  \cite{An74}, p.115).

\

On the other hand, 
for $n=2$ and for two cusp eigenforms  $f$ and $g$ for $\Sp_2(\Z)$  of weights $k,l$, $k>l+1$, 
$\Lambda(Sp(f)\otimes Sp(g), s)=\ep(f,g)\Lambda(Sp(f)\otimes Sp(g), 2k+2l-5-s)$, $|\ep(f,g)|=1$, where
\begin{align*}&
\Lambda(Sp(f)\otimes Sp(g), s)=\Gamma_\CC(s)\Gamma_\CC(s-l+2)\Gamma_\CC(s-l+1)\Gamma_\CC(s-k+2)
\\ & \times
\Gamma_\CC(s-k+1)\Gamma_\CC(s-2l+3)\Gamma_\CC(s-k-l+2)\Gamma_\CC(s-k-l+3)
\\ & \times
L(Sp(f)\otimes Sp(g), s).
\end{align*}
We used here the Gauss duplication formula $\Gamma_\CC(s)=\Gamma_\RR(s)\Gamma_\RR(s+1)$. 
Notice that  $a_+=a_-=1$ in this case,
and the conjectural motive $M(Sp(f))\otimes M(Sp(g))$ does not admit critical values.

\section{A holomorphic lifting from $GSp_2 \times GSp_2$ to $GSp_4$: 
a conjecture} 
(compare with constructions in \cite{BFG06}, \cite{BFG92}, \cite{Jia96} for generic automorphic forms).

Our computation makes it possible to compare the 
spinor Hecke series of genus 4 computed in \cite{VaSp4}
(in variables $u_0, u_1, u_2, u_3, u_4$)
with the Rankin product of two Hecke series of genus 2 (in variables $x_0, x_1, x_2, y_0, y_1, y_2$).
It follows from our computation that if we make the substitution $u_0=x_0y_0, u_1=x_1, u_2=x_2, u_3=y_1, u_4=y_2$
then the denominator of the series
$$
\sum _{\delta =0}^{\infty } 
\,{\Omega^{(4)} _{u}}(\mathbf{T}(p^{\delta }))X^{\delta } 
$$
coincides with the denominator of the Rankin product
$$
\sum _{\delta =0}^{\infty } 
\,{\Omega^{(2)} _{x}}(\mathbf{T}(p^{\delta }))\cdot \Omega^{(2)} _{y}(\mathbf{T}(p^{\delta }))\,X^{\delta }
\in\Q[x_0, x_1, x_2, y_0, y_1, y_2][\![X]\!].
$$
On the basis of this equality we would like to push forward the following

\begin{conj}[on a lifting from $GSp_2 \times GSp_2$ to $GSp_4$ (of genus four)]\label{lift}
Let $f$ and $g$ be two Siegel modular forms of genus 2 and of weights $k> 4$ and $l=k-2$. 
Then there exists a Siegel modular form $F$ of genus 4 and of weight $k$ with the Satake parameters
$$\gamma_0=\alpha_0\beta_0, \gamma_1=\alpha_1, \gamma_2=\alpha_2, \gamma_3=\beta_1, \gamma_4=\beta_2,
$$ 
for a suitable choice of  Satake's parameters 
$\alpha_0, \alpha_1, \alpha_2$ and  $\beta_0, \beta_1, \beta_2$ of $f$ and $g$.

\end{conj}
\begin{rema}\label{Ike}
An evidence for the conjecture comes from Ikeda-Miyawaki constructions (\cite{Ike01}, \cite{Mur02}, \cite{Ike06}): let $k$ be an even positive integer, 
${h}\in S_{2k}(\Gamma_1)$ a normalized Hecke eigenform of weight $2k$, 
$F_2(h)\in S_{k+1}(\Gamma_2) =Maass({h})$ the Maass lift of ${h}$,  and
$F_{2n}\in S_{k+n}(\Gamma_{2n})$ the Ikeda lift of ${h}$
(we assume $k\equiv n\bmod 2$, $n\in \N$).


Next let
$f\in S_{k+n+r}(\Gamma_{r})$  be an arbitrary   Siegel cusp eigenform of genus $r$ and weight $k+n+r$,  with $n,r\ge 1$.
Then according to Ikeda-Miyawaki (see \cite{Ike06}) there exists a Siegel eigenform
${\cal F}_{{h},f}\in S_{k+n+r}(\Gamma_{2n+r})$
such that 
\begin{align}\label{IkeMi}
L(s,{\cal F}_{{h},f}, St) = L(s,f,St)
\prod_{j=1}^{2n}L(s+k+n-j,{h}) 
\end{align}
(under a non-vanishing condition, see Theorem 2.3 at p.63 in \cite{Mur02}).
The form ${\cal F}_{{h},f}$ is given by the integral 
$$
{\cal F}_{{h}, f}(Z)=\langle F_{2n+2r}(\diag(Z,Z'), f(Z')\rangle_{Z'}\,.
$$
If we take $n=1,r=2$, $k:=k+1$ then an example of the validity of the conjecture is given by $g=F_2(h)$, 
$$
(f,g)=(f, F_2(h))\mapsto {\cal F}_{f, {h}}\in S_{k+3}(\Gamma_{4}), \ (f, g)=(f, F_2(h))\in  
S_{k+3}(\Gamma_{2})\times S_{k+1}(\Gamma_{2}).
$$
\end{rema}
\begin{rema}
Notice that the Satake parameters of the Ikeda lift $F=F_{2m}(h)$ of $h$ can be taken in the form 
$\beta_0, \beta_1, \cdots, \beta_{2m}$, where 
$$
\beta_0=p^{{m}k-{m}({m}+1)/2}, 
\beta_i=\alpha p^{i-1/2},\ \  \beta_{{m}+i}=\alpha^{-1} p^{i-1/2}, \ (i=1, \cdots,{m})
$$
and 
$$
(1-\alpha p^{k-1/2}X)(1-\alpha^{-1}p^{k-1/2}X)=1-a(p)X+p^{2k-1}X^2, h=\sum_{n=1}^\infty a(n)q^n 
$$
see \cite{Mur02}.
\end{rema}


The $L$-function of degree 16 in Conjecture \ref{lift}
is related to the tensor product $L$-function in \cite{Jia96}.
In the example of Remark \ref{Ike}
it coincides with the product of two shifted $L$-functions of degree 8 of Boecherer-Heim 
\cite{BoeH06}.


\begin{conj}[on a lifting from $GSp_{2m} \times GSp_{2m}$ to $GSp_{4m}$]\label{lift2m}

 Here is a version of Conjectire \ref{lift} for any even genus $r=2m$.
Let $f$ and $g$ be two Siegel modular forms of genus $2m$ and of weights $k> 2m$ and $l=k-2m$.
Then there exists a Siegel modular form $F$ of genus $4m$ and of weight $k$ with the Satake parameters
$\gamma_0=\alpha_0\beta_0$, $\gamma_1=\alpha_1$, $\gamma_2=\alpha_2, \cdots, \gamma_{2m}=\alpha_{2m}$, 
$\gamma_{2m+1}=\beta_1, \cdots, \gamma_{4m}=\beta_{2m}$ for suitable choices
$\alpha_0$, $\alpha_1, \cdots, \alpha_{2m}$ and  $\beta_0$, $\beta_1, \cdots, \beta_{2m}$
of Satake's parameters of $f$ and $g$.


One readily checks that the Hodge types of $M(Sp(f))\otimes M(Sp(g))$  and $M(Sp(F))$ are again the same (of rank $2^{4m}$)
(it follows from the above description (\ref{Hodge}), and from Künneth's-type formulas). 

\end{conj}

An evidence for this version of the conjecture comes again from Ikeda-Miyawaki constructions (\cite{Ike01}, \cite{Mur02}, \cite{Ike06}): 
let $k$ be an even positive integer, 
${h}\in S_{2k}(\Gamma_1)$ a normalized Hecke eigenform of weight $2k$, 
$F_{2n}\in S_{k+n}(\Gamma_{2n})$ the Ikeda lift of ${h}$ of genus $2n$
(we assume $k\equiv n\bmod 2$, $n\in \N$).

Next let
$f\in S_{k+n+r}(\Gamma_{r})$  be an arbitrary   Siegel cusp eigenform of genus $r$ and weight $k+n+r$,  with $n,r\ge 1$.
If we take in (\ref{IkeMi}) $n=m, r=2m$, $k:=k+m$, $k+n+r:=k+3m$, 
then an example of the validity of this version of the conjecture is given by 
$$
(f, g)=(f, F_{2m}(h))\mapsto {\cal F}_{{h},f}\in S_{k+3m}(\Gamma_{4m}), \ 
(f, g)=(f, F_{2m})\in   S_{k+3m}(\Gamma_{2m})\times S_{k+m}(\Gamma_{2m}).
$$

\

Another evidence comes from Siegel-Eiseinstein series 
$$
f=E^{2m}_k  \mbox{ and  }g=E^{2m}_{k-2m}
$$
 of  even genus $2m$ and weights $k$ and $k-2m$:
we have then
\begin{align*}&
\alpha_0=1, \alpha_1=p^{k-2m}, \cdots, \alpha_{2m}=p^{k-1},
\\ &
\beta_0=1, \beta_1=p^{k-4m}, \cdots, \beta_{2m}=p^{k-2m-1},
\end{align*}
then we have that
$$
\gamma_0=1, \gamma_1=p^{k-4m}, \cdots, \gamma_{2m}=p^{k-1},
$$
are the Satake parameters of the Siegel-Eisenstein series $F=E^{4m}_{k}$. 

\appendix\label{coeff}
\section{Appendix: Coefficients of the polynomials $R(X)$ and $S(X)$}
\begin{figure}[h]
\caption{\label{NP} Newton polygons of $R(X)$ and $S(X)$ with respect to powers of $p$ and $X$,  of heights 34 and 48, resp.} 
\begin{center}
\includegraphics[width=5cm,height=5cm]{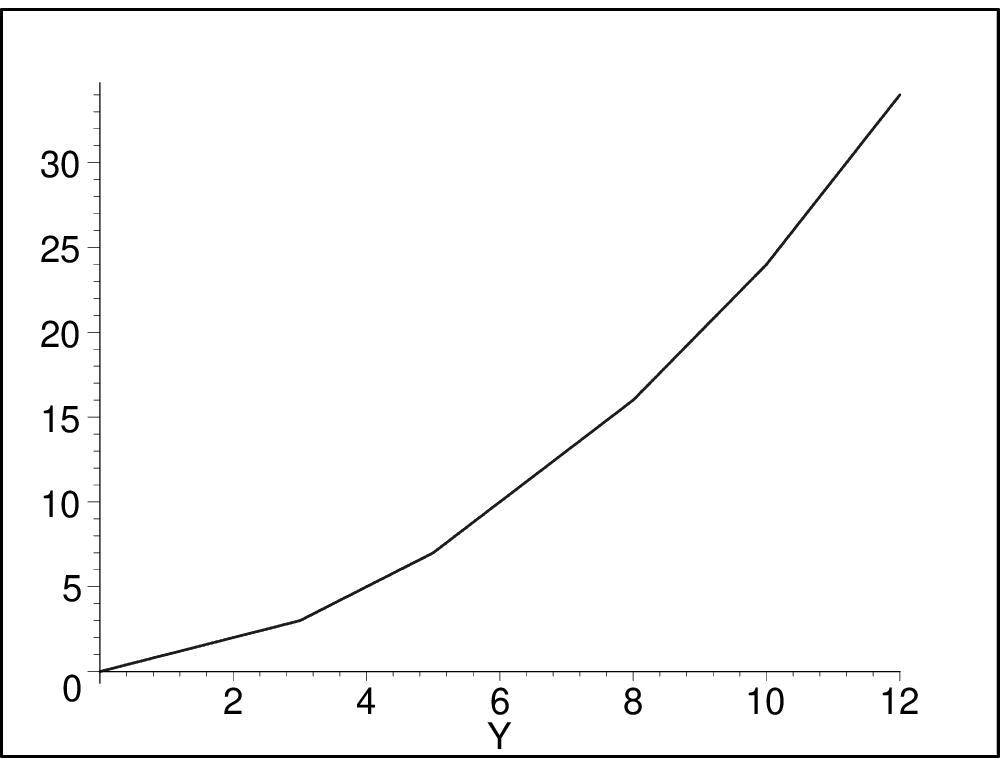}\hskip0.5cm
\includegraphics[width=5cm,height=5cm]{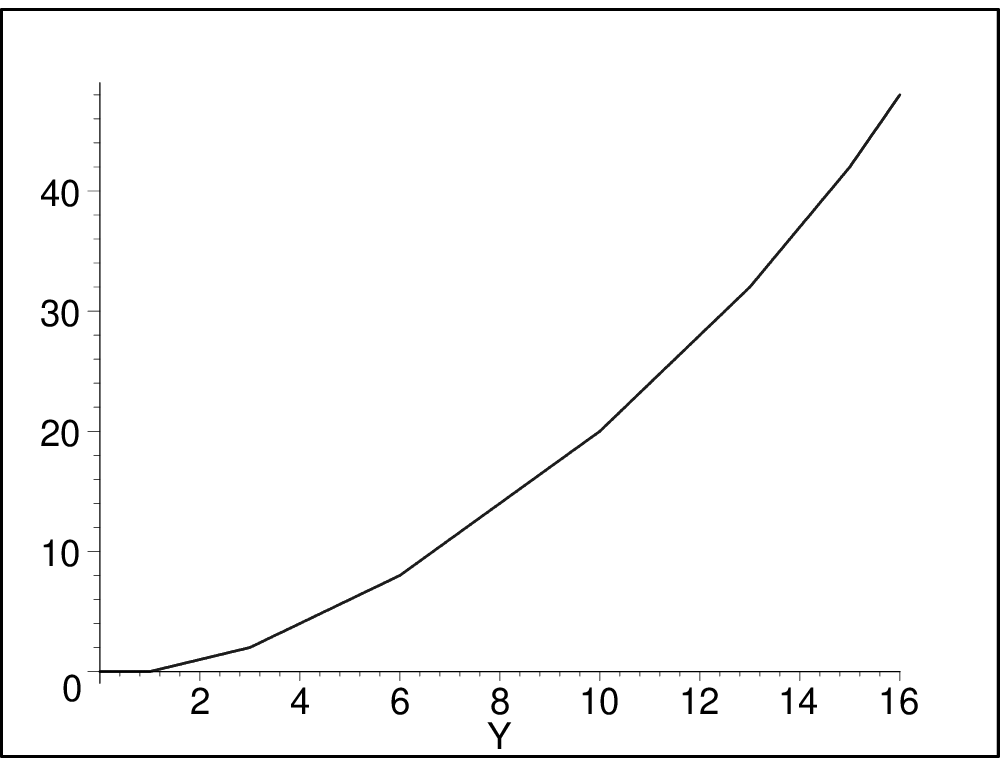}
\end{center}
\end{figure}

We give here explicit expressions for the coefficients of the polynomials $R(X)$ and $S(X)$
from Theorem \ref{ThRS}. From these formulas one can observe some nice divisibility properties (by certain powers of $p$ and 
the elements $\Tc\otimes \Tc\in{\mathcal L}_{2,\Z}\otimes{\mathcal L}_{2,\Z}$):
\begin{align*}
R(X)&=1+r_2X^2+\cdots +r_{10}X^{10}+r_{12}X^{12} \in{\mathcal L}_{2,\Z}\otimes {\mathcal L}_{2,\Z}[X] \mbox{ with } r_{1}=r_{11}=0,
\\ 
S(X)
&=1+s_1X+\cdots +s_{16}X^{16}
\\ &\nonumber
=
1 - (\Ta \otimes \Ta) X+\cdots +(p^6\Tc\otimes \Tc)^{8}X^{16}\in{\mathcal L}_{2,\Z}\otimes{\mathcal L}_{2,\Z}[X], 
\end{align*}
 with $r_i$ and $s_i$  given as follows
\begin{align*}&
r_2=
p^2((2p-1)(p^2+1)
\Tc\otimes \Tc 
-(p^2-p+1)
(
\Tb\otimes \Tc 
+\Tc\otimes \Tb) \\ & \nonumber \quad
-(                
\Tb\otimes \Tb 
+               
\Ta^2\otimes \Tc 
+\Tc\otimes \Ta^2  
),
\\ &
\\ &
r_3=
p^3(p+1)(2\Tc\otimes\Tc+\Tb\otimes\Tc+\Tc\otimes\Tb)\Ta\otimes\Ta \, ,
\\ &
\\ &
r_4=
%
-p^5(
(p^7+2p^6-2p^5+6p^4+p^3+6p^2+p+2)  \Tc^2\otimes\Tc^2                             \\ & \nonumber\quad
-(p^2+1)(p^3-3p^2-p-3)            (\Tb\otimes\Tc+\Tc\otimes\Tb)\Tc\otimes\Tc     \\ & \nonumber\quad
+(p+4)(p^2+1)                      \Tb\Tc\otimes\Tb\Tc
-(p^3-p^2-1)                      (\Tb^2\otimes\Tc^2+\Tc^2\otimes\Tb^2)          \\ & \nonumber\quad
+                                 (\Tb\otimes\Tc+\Tc\otimes\Tb)\Tb\otimes\Tb
-p(p^3+2p^2-p+2)                 (\Ta^2\otimes\Tc\\ & \nonumber \quad +\Tc\otimes\Ta^2)\Tc\otimes\Tc
-2p                               (\Ta^2\otimes\Tb+\Tb\otimes\Ta^2)\Tc\otimes\Tc \\ & \nonumber\quad
+p^2                              (\Ta^2\Tb\otimes\Tc^2+\Tc^2\otimes\Ta^2\Tb)
+(p+2)                             \Ta^2\Tc\otimes\Ta^2\Tc
) \, ,
\end{align*}
\begin{align*} 
&
r_5=
-p^7(
2(p+1)(2p^4-p^3+p^2-1) \Tc\otimes\Tc
+(p+1)(p-2)           (\Tb\otimes\Tc+\Tc\otimes\Tb)       \\ & \nonumber\quad
-2                     \Tb\otimes\Tb              
-p(p+1)               (\Ta^2\otimes\Tc+\Tc\otimes\Ta^2)
) \Ta\Tc\otimes\Ta\Tc \, ,
\\ & 
\\ &
r_6=
-p^{10}\,(
 p\,(p^2+1)(p^5-2p^3-8p^2-p-4) \Tc^3\otimes\Tc^3                                       \\ & \quad 
-p\,(p^5+4p^4+2p^3+12p^2+p+6) (\Tb\otimes\Tc+\Tc\otimes\Tb)\Tc^2\otimes\Tc^2           \\ & \quad 
+p\,(p-4)(p^2+1)               \Tb\Tc^2\otimes\Tb\Tc^2                                 \\ & \quad 
-p\,(p+4)(p^2+1)              (\Tb^2\otimes\Tc^2+\Tc^2\otimes\Tb^2)\Tc\otimes\Tc       \\ & \quad 
-p\,                          (\Tb\otimes\Tc+\Tc\otimes\Tb)\Tb\Tc\otimes\Tb\Tc         \\ & \quad 
-p\,                          (\Tb^3\otimes\Tc^3+\Tc^3\otimes\Tb^3)                    \\ & \quad 
-(p^5-4p^2-p-2)               (\Ta^2\otimes\Tc+\Tc\otimes\Ta^2)\Tc^2\otimes\Tc^2       \\ & \quad 
+(p^2+3)                      (\Ta^2\otimes\Tb+\Tb\otimes\Ta^2)\Tc^2\otimes\Tc^2       \\ & \quad 
+                             (\Ta^2\Tc\otimes\Tb^2+\Tb^2\otimes\Ta^2\Tc)\Tc\otimes\Tc \\ & \quad 
+(p^3+3p^2+p+1)               (\Ta^2\Tb\otimes\Tc^2+\Tc^2\otimes\Ta^2\Tb)\Tc\otimes\Tc \\ & \quad 
+                             (\Ta^2\otimes\Tc+\Tc\otimes\Ta^2)\Tb\Tc\otimes\Tb\Tc     \\ & \quad 
+(p^2+1)                       \Ta^2\Tc^2\otimes\Ta^2\Tc^2) \, ,
\\ & 
\\ & 
r_7=
-p^{13}\, (
 2(p+1)(p^3+p-1) \Tc\otimes\Tc
-(p+1)(p^2-2p+2)(\Tb\otimes\Tc+\Tc\otimes\Tb)     \\ & \quad
- 2              \Tb\otimes\Tb
-(p+1)          (\Ta^2\otimes\Tc+\Tc\otimes\Ta^2)
) \Ta\Tc^2\otimes\Ta\Tc^2 \, ,
\\ & 
\\ & 
r_8=
-p^{16}(
 p\,(2p^6+3p^5+6p^4-p^3+6p^2-p+2) \Tc^2\otimes\Tc^2                             \\ & \quad
+p\,(p^2+1)(p^3+3p^2-p+3)        (\Tb\otimes\Tc+\Tc\otimes\Tb)\Tc\otimes\Tc     \\ & \quad
+p\,(p+4)(p^2+1)                  \Tb\Tc\otimes\Tb\Tc                           \\ & \quad
+p\,(p^2-p+1)                    (\Tb^2\otimes\Tc^2+\Tc^2\otimes\Tb^2)          \\ & \quad
+p\,                             (\Tb\otimes\Tc+\Tc\otimes\Tb)\Tb\otimes\Tb     \\ & \quad
-p\,(2p^3+p^2+2p-1)              (\Ta^2\otimes\Tc+\Tc\otimes\Ta^2)\Tc\otimes\Tc \\ & \quad
-2p^2                            (\Ta^2\otimes\Tb+\Tb\otimes\Ta^2)\Tc\otimes\Tc \\ & \quad
+p\,                             (\Ta^2\Tb\otimes\Tc^2+\Tc^2\otimes\Ta^2\Tb)    \\ & \quad
+(2p+1)                           \Ta^2\Tc\otimes\Ta^2\Tc)\Tc^2\otimes\Tc^2 \, ,
\\ & 
\\ &
r_9=
p^{20}(p+1)(2\Tc\otimes\Tc+\Tb\otimes\Tc+\Tc\otimes\Tb)\Ta\Tc^3\otimes\Ta\Tc^3
\\ & 
\\ & 
r_{10}=
p^{24}(
 (p^2+1)(p^4+2p^3-p^2-1) \Tc\otimes\Tc
+(p^3-p^2-1)            (\Tb\otimes\Tc+\Tc\otimes\Tb)  \\ & \quad
-                        \Tb\otimes\Tb
-p^2                    (\Ta^2\otimes\Tc+\Tc\otimes\Ta^2) ) \Tc^4\otimes\Tc^4\,,
\\ & 
r_{11}=0,
\\ & 
\\ & 
r_{12}=
p^{34} 
\Tc^6 \otimes \Tc^6\,.
\end{align*}
As for the coefficients of $S(X)$, one has
\begin{align*}
S(X)=
1 - (\Ta \otimes \Ta) X+\cdots +(p^6\Tc\otimes \Tc)^{8}X^{16}\in{\mathcal L}_{2,\Z}\otimes{\mathcal L}_{2,\Z}[X], 
\end{align*}
where
\begin{align*}& 
s_1= -\Ta\otimes\Ta \,,
\\ & 
\\ &
s_2=
- p(2\,p\,(p^{2} + 1)^{2}\,
\Tc\otimes\Tc
+ 2\,p\,(p^{2} + 1)\,
(\Tb\otimes\Tc+\Tc\otimes\Tb)
\\ &\quad 
\mbox{} + 2\,p\,
\Tb\otimes\Tb
 - (p^{2} + 1)\,
(\Ta^2\otimes\Tc+\Tc\otimes\Ta^2)
%
\\   & \quad
- 
(\Ta^2\otimes\Tb+\Tb\otimes\Ta^2)) \,,
\\ & 
\\ & 
s_3= p^{2} 
((2\,p^{4} + 4\,p^{2} - 1)\,
\Tc\otimes\Tc
+ (2\,p^{2} - 1)\,
(\Tb\otimes\Tc+\Tc\otimes\Tb)
\\ &\quad 
-\Tb\otimes\Tb
 - p\,
(\Ta^2\otimes\Tc+\Tc\otimes\Ta^2)
) 
\Ta\otimes\Ta
) \,,
\\ & 
\\ & 
s_{4}=
p^{4}((p^{8} + 12\,p^{6} + 10\,p^{4} + 4\,p^{2} + 1)\,
\Tc^2\otimes\Tc^2
\\ &\quad 
\mbox{} + 2\,(3\,p^{6} + 5\,p^{4} + 3\,p^{2} + 1)\,
(\Tb\otimes\Tc+\Tc\otimes\Tb)
\Tc\otimes\Tc
\\ &\quad 
\mbox{} + 4\,(p^{2} + 1)^{2}\,
\Tb\Tc\otimes\Tb\Tc
\\ &\quad
+ 
(3\,p^{4} + 2\,p^{2} + 1)\,
(\Tb^2\otimes\Tc^2+\Tc^2\otimes\Tb^2)
\\ &\quad 
\mbox{} + 2\,(p^{2} + 1)\,
(\Tb\otimes\Tc+\Tc\otimes\Tb)
\Tb\otimes\Tb
\\ &\quad
+ 
\Tb^2\otimes\Tb^2
\\ &\quad 
\mbox{} - 2\,p\,(p^{4} + 4\,p^{2} + 1)\,
(\Ta^2\otimes\Tc+\Tc\otimes\Ta^2)
\,
\Tc\otimes\Tc
\\ &\quad 
- 4\,p\,(p^{2} + 1)\,
(\Ta^2\otimes\Tb+\Tb\otimes\Ta^2)
\Tc\otimes\Tc
%
\\ &\quad 
- 2\,p\,
(\Ta^2\Tc\otimes\Tb^2+\Tb^2\otimes\Ta^2\Tc)
\\ &\quad
- 4\,p^{3}\,
(\Ta^2\Tb\otimes\Tc^2+\Tc^2\otimes\Ta^2\Tb)
\\ &\quad 
\mbox{} + (p^{2} + 2)\,
\Ta^2\Tc\otimes\Ta^2\Tc
\\ &\quad
+ 
(\Tb\otimes\Tc+\Tc\otimes\Tb)
\Ta^2\otimes\Ta^2
\\ &\quad 
+ p^{2}\,
(\Ta^4\otimes\Tc^2+\Tc^2\otimes\Ta^4)
) \,,
\\ & 
\\ & 
s_{5}= - p^{6}( 
(6\,p^{6} + 2\,p^{4} - p^{2} + 2)\,
\Tc^2\otimes\Tc^2
\\ &\quad 
+ (p^{4} - p^{2} + 3)\,
(\Tb\otimes\Tc+\Tc\otimes\Tb)
\Tc\otimes\Tc
\\ &\quad 
\mbox{} + (3\,p^{2} + 4)\,
\Tb\Tc\otimes\Tb\Tc
\\ &\quad 
- (2\,p^{2} - 1)\,
(\Tb^2\otimes\Tc^2+\Tc^2\otimes\Tb^2)
\\ &\quad 
+ (
\Tb\otimes\Tc+\Tc\otimes\Tb
)\,
\Tb\otimes\Tb
\\ &\quad 
\mbox{} - p\,(2\,p^{2} + 1)\,
(\Ta^2\otimes\Tc+\Tc\otimes\Ta^2)
\Tc\otimes\Tc
\\ &\quad 
- 2\,p\,
(\Ta^2\otimes\Tb+\Tb\otimes\Ta^2) 
\Tc\otimes\Tc
\\ &\quad 
+ p\,
(\Ta^2\Tb\otimes\Tc^2+\Tc^2\otimes\Ta^2\Tb)
\\ &\quad 
+ 
\Ta^2\Tc\otimes\Ta^2\Tc
)
\Ta\otimes\Ta
) \,,
\end{align*}

\begin{align*}&  
s_{6}= - p^{8}(2\,p^{2}\,(p^{8} + 6\,p^{6} + 
11\,p^{4} + 8\,p^{2} + 2)\,
\Tc^3\otimes\Tc^3
\\ &\quad 
\mbox{} + 2\,p^{2}\,(5\,p^{4} + 12\,p^{2} + 6)\,
\Tb\Tc^2\otimes\Tb\Tc^2
\\ &\quad 
\mbox{} + (3\,p^{4} + 10\,p^{2} - 1)\,
\Ta^2\Tc^2\otimes\Ta^2\Tc^2
- 
\Ta^2\Tb\Tc\otimes\Ta^2\Tb\Tc
\\ &\quad 
\mbox{} + 2\,p^{2}\,(3\,p^{6} + 11\,p^{4} + 12\,p^{2} + 4)\,
(
\Tb\otimes\Tc+\Tc\otimes\Tb
)\,
\Tc^2\otimes\Tc^2
\\ &\quad 
\mbox{} + 6\,p^{2}\,(p^{2} + 1)^{2}\,
(\Tb^2\otimes\Tc^2+\Tc^2\otimes\Tb^2)
\Tc\otimes\Tc
\\ &\quad
+ 6\,p^{2}\,(p^{2} + 1)\,(
\Tb\otimes\Tc+\Tc\otimes\Tb
)\,
\Tb\Tc\otimes\Tb\Tc
\\ &\quad 
\mbox{} + 2\,p^{2}\,(p^{2} + 1)\,
(\Tb^3\otimes\Tc^3+\Tc^3\otimes\Tb^3)
\\ &\quad
+ 2\,p^{2}\,
(\Tb^2\otimes\Tc^2+\Tc^2\otimes\Tb^2)
\Tb\otimes\Tb
\\ &\quad
- p\,(5\,p^{6} + 13\,p^{4} + 10\,p^{2} + 2)\,
(\Ta^2\otimes\Tc+\Tc\otimes\Ta^2)
\Tc^2\otimes\Tc^2
\\ &\quad 
\mbox{} - p\,(7\,p^{4} + 12\,p^{2} + 4)\,
(\Ta^2\otimes\Tb+\Tb\otimes\Ta^2)
\Tc^2\otimes\Tc^2
\\ &\quad
- 3p\,(\,p^{2} + 1)\,
(\Ta^2\Tc\otimes\Tb^2+\Tb^2\otimes\Ta^2\Tc)
\Tc\otimes\Tc
\\ & \quad 
\mbox{} - p\,(
\Ta^2\Tc^2\otimes\Tb^3+\Tb^3\otimes\Ta^2\Tc^2
) 
\\ &\quad
- 2\,p
\,(3\,p^{4} + 4\,p^{2} + 1)\,
(\Ta^2\Tb\otimes\Tc^2+\Tc^2\otimes\Ta^2\Tb)
\Tc\otimes\Tc
\\ &\quad 
\mbox{} - 2\,p\,(3\,p^{2} + 1)\,
(\Ta^2\otimes\Tc+\Tc\otimes\Ta^2)
\Tb\Tc\otimes\Tb\Tc
\\ &\quad
 - p\,(p^{2} + 1)\,
(\Ta^2\Tb^2\otimes\Tc^3+\Tc^3\otimes\Ta^2\Tb^2)
\\ &\quad
 - p\,
(\Ta^2\Tb\otimes\Tc^2+\Tc^2\otimes\Ta^2\Tb)
\Tb\otimes\Tb
\\ &\quad 
\mbox{} + (5\,p^{2} - 1)\,(
\Tb\otimes\Tc+\Tc\otimes\Tb
)\,
\Ta^2\Tc\otimes\Ta^2\Tc
\\ &\quad
+ 2\,p^{2}\,(p^{2} + 1)\,
(\Ta^4\otimes\Tc^2+\Tc^2\otimes\Ta^4)
\Tc\otimes\Tc
\\ &\quad 
\mbox{} 
+ 2\,p^{2}\,(
\Ta^4\otimes\Tb\Tc+\Tb\Tc\otimes\Ta^4
)\,
\Tc\otimes\Tc
\\ &\quad 
\mbox{} - p\,(
\Ta^4\otimes\Ta^4\Tc+\Ta^4\Tc\otimes\Ta^4
)\,
\Tc\otimes\Tc 
) \,,
\\ & 
\\ & 
  s_{7}= p^{11}( 
p\,(5\,p^{6} - 2\,p^{4} + 2)\,
\Ta\Tc^3\otimes\Ta\Tc^3
 \\ &\quad 
+ 8\,p\,
\Ta\Tb\Tc^2\otimes\Ta\Tb\Tc^2
 \\ &\quad 
\mbox{} + p\,
\Ta^3\Tc^2\otimes\Ta^3\Tc^2
 \\ &\quad 
- p\,(p^{4} - 3)\,(
\Tb\otimes\Tc+\Tc\otimes\Tb
)\,
\Ta\Tc^2\otimes\Ta\Tc^2
\\ &\quad 
 - p\,
(\Tb^2\otimes\Tc^2+\Tc^2\otimes\Tb^2)
\Ta\Tc\otimes\Ta\Tc
\\ &\quad 
+ 2\,p\,(
\Tb\otimes\Tc+\Tc\otimes\Tb
)\,
\Ta\Tb\Tc\otimes\Ta\Tb\Tc
\\ &\quad 
\mbox{} - p\,
(\Tb^3\otimes\Tc^3+\Tc^3\otimes\Tb^3)
\Ta\otimes\Ta
\\ &\quad 
- (3\,p^{4}
 - 3\,p^{2} + 2)\,
(\Ta^2\otimes\Tc+\Tc\otimes\Ta^2)
\Ta\Tc^2\otimes\Ta\Tc^2
\\ &\quad 
\mbox{} + (p^{2} - 3)\,
(\Ta^2\otimes\Tb+\Tb\otimes\Ta^2)
\Ta\Tc^2\otimes\Ta\Tc^2
\\ &\quad 
- 
(\Ta^2\Tc\otimes\Tb^2+\Tb^2\otimes\Ta^2\Tc)
\Ta\Tc\otimes\Ta\Tc
\\ &\quad 
\mbox{} 
+ (2\,p^{2} - 1)\,
(\Ta^2\Tb\otimes\Tc^2+\Tc^2\otimes\Ta^2\Tb)
\Ta\Tc\otimes\Ta\Tc
\\ &\quad 
-(\Ta^2\otimes\Tc+\Tc\otimes\Ta^2)
\Ta\Tb\Tc\otimes\Ta\Tb\Tc
) \,,
\end{align*}

\begin{align*} &  s_{8}= p^{14}( 
2\,p^{2}\,(2\,p^{8} + 4\,p^{6} + 14\,p^{4} + 12\,p^{2} + 3)\,
\Tc^4\otimes\Tc^4
\\ &\quad 
\mbox{} + 4\,p^{2}\,(p^{6} + 7\,p^{4} + 9\,p^{2} + 3)\,(
\Tb\otimes\Tc+\Tc\otimes\Tb
)\,
\Tc^3\otimes\Tc^3
\\ &\quad 
\mbox{} + 16\,p^{2}\,(p^{2} + 1)^{2}\,
\Tb\Tc^3\otimes\Tb\Tc^3
\\ &\quad 
+ 2\,p^{2}\,(3\,p^{4} + 10\,
p^{2} + 5)\,
(\Tb^2\otimes\Tc^2+\Tc^2\otimes\Tb^2)
\,
\Tc^2\otimes\Tc^2
\\ &\quad 
\mbox{} + 8\,p^{2}\,(p^{2} + 1)\,(
\Tb\otimes\Tc+\Tc\otimes\Tb
)\,
\Tb\Tc^2\otimes\Tb\Tc^2
\\ &\quad 
\mbox{} + 4\,p^{2}\,
\Tb^2\Tc^2\otimes\Tb^2\Tc^2
\\ &\quad 
+ 4\,p^{2}\,(p^{2} + 1)\,
(\Tb^3\otimes\Tc^3+\Tc^3\otimes\Tb^3)
\,
\Tc\otimes\Tc
\\ &\quad 
\mbox{} + p^{2}\,(
\Tb^4\otimes\Tc^4+\Tc^4\otimes\Tb^4
) 
\\ &\quad 
- 4\,p\,(2\,p^{6} + 3\,p^{4} + 4\,p^{2} + 1)\,
(\Ta^2\otimes\Tc+\Tc\otimes\Ta^2)
\Tc^3\otimes\Tc^3
\\ &\quad 
\mbox{} - 8\,p\,(p^{2} + 1)^{2}\,
(\Ta^2\otimes\Tb+\Tb\otimes\Ta^2)
\Tc^3\otimes\Tc^3
\\ &\quad 
 - 4\,p\,(p^{2} + 1)\,
(\Ta^2\Tc\otimes\Tb^2+\Tb^2\otimes\Ta^2\Tc)
\Tc^2\otimes\Tc^2
\\ &\quad 
\mbox{} - 4\,p\,(p^{4} + 4\,p^{2} + 1)\,
(\Ta^2\Tb\otimes\Tc^2+\Tc^2\otimes\Ta^2\Tb)
\Tc^2\otimes\Tc^2
\\ &\quad 
 - 8\,p\,(p^{2} + 1)\,
(\Ta^2\otimes\Tc+\Tc\otimes\Ta^2)
\Tb\Tc^2\otimes\Tb\Tc^2
 \\ &\quad 
\mbox{} - 4\,p\,
(\Ta^2\otimes\Tb+\Tb\otimes\Ta^2)
\Tb\Tc^2\otimes\Tb\Tc^2
\\ &\quad 
- 4\,p^{3}\,
(\Ta^2\Tb^2\otimes\Tc^3+\Tc^3\otimes\Ta^2\Tb^2) + 
\Tc\otimes\Tc 
\\ &\quad 
\mbox{} + 2\,(5\,p^{4} + 2\,p^{2} + 2)\,
\Ta^2\Tc^3\otimes\Ta^2\Tc^3
\\ &\quad 
\mbox{} + 2\,(p^{2} + 2)\,(
\Tb\otimes\Tc+\Tc\otimes\Tb
)\,
\Ta^2\Tc^2\otimes\Ta^2\Tc^2
\\ &\quad 
\mbox{} + 2\,
\Ta^2\Tb\Tc^2\otimes\Ta^2\Tb\Tc^2
\\ &\quad 
+ (\Tb^2\otimes\Tc^2+\Tc^2\otimes\Tb^2)
\Ta^2\Tc\otimes\Ta^2\Tc
 \\ &\quad 
\mbox{} + (3\,p^{4} + 2\,p^{2} + 1)\,
(\Ta^4\otimes\Tc^2+\Tc^2\otimes\Ta^4)
\Tc^2\otimes\Tc^2
\\ &\quad 
\mbox{} + 2\,(p^{2} + 1)\,(
\Ta^4\otimes\Tb\Tc+\Tb\Tc\otimes\Ta^4
)\,
\Tc^2\otimes\Tc^2
\\ &\quad 
\mbox{} + (
\Ta^4\otimes\Tb^2+\Tb^2\otimes\Ta^4
)\,
\Tc^2\otimes\Tc^2
\\ &\quad 
 - 2\,p\,
(\Ta^2\otimes\Tc+\Tc\otimes\Ta^2)
\Ta^2\Tc^2\otimes\Ta^2\Tc^2
) \,.
\end{align*}
Then we find the remaining coefficients $s_9, \cdots, s_{16}$,
using an easy  functional equation   (similar to \cite{An87}, p.164, (3.3.79)):   
$$
s_{16-i}=(p^6\Tc\otimes \Tc)^{8-i}s_i \ \ \ (i=0, \cdots, 8).
$$
To conclude with, we give the Newton polygons of $R(X)$ and $S(X)$ with respect to powers of $p$ and $X$
(see Figure \ref{NP}). It follows from our computation that all slopes are {\it integral}. 
We hope that these polygons could  help to find some geometric objects 
attached to
the polynomials  $R(X)$ and $S(X)$, 
in the spirit of a recent work of C.~Faber and G.~Van~Der~Geer, \cite{FVdG}.
\subsection*{Acknowledgement}
We are very grateful to the Institute Fourier
 (UJF, Grenoble-1) 
for  the  permanent excellent working environment. 

\

It is a great pleasure for us to thank  Siegfried Boecherer, Yasutaka Ihara, Solomon Friedberg and Mikhail Tsfasman
 for valuable  discussions and observations. 

\

Our special thanks go to Yuri Ivanovich Manin, 
 for providing us with  advice and encouragement.

\bibliographystyle{plain}

\end{document}